\documentclass[12pt]{amsart}
\usepackage{amsmath,amsthm,amsfonts,amssymb}
\usepackage[all]{xy}
\usepackage{amscd}

\DeclareMathOperator{\ad}{ad}

\DeclareMathOperator{\Gr}{Gr}
\DeclareMathOperator{\Hom}{Hom}
\DeclareMathOperator{\Ht}{ht}
\DeclareMathOperator{\im}{im}
\DeclareMathOperator{\Ker}{ker}
\DeclareMathOperator{\Span}{Span}

\newcommand{\coh}{{\rm H}}

\newcommand{\HH}{{\rm H}_{h}}

\newcommand{\HB}{{\rm B}_{h}}
\newcommand{\HZ}{{\rm Z}_{h}}
\newcommand{\Hb}{\widehat{\rm H}_{b}}
\newcommand{\Bb}{\widehat{\rm B}_{b}}
\newcommand{\Zb}{\widehat{\rm Z}_{b}}

\def\a{\mathbf{a}}
\def\b{\mathbf{b}}
\def\x{\mathbf{x}}
\def\y{\mathbf{y}}
\def\ot{\otimes}
\def\op{{\operatorname{op}}}
\def\cop{{\operatorname{cop}}}
\def\p #1#2{\langle #1|#2\rangle}
\def\l #1#2{{^#1 #2}}
\def\Hom{\operatorname{Hom}}
\def\Homk{\operatorname{Hom}_k}
\def\de{\delta}
\def\De{\Delta}

\def\sx{\sigma^X}
\def\sb{\sigma^B}
\def\na{\partial}
\def\nax{\na^X}
\def\nab{\na^B}
\def\nah{\na^h}
\def\nac{\na^c}

\def\ep{\varepsilon}

\def\hy{{\widehat\y}}
\def\hb{{\widehat\b}}
\def\C{\mathbf{C}}
\def\B{\mathbf{B}}
\def\A{\mathbf{A}}

\def\res{\operatorname{res}}
\def\Tot{\operatorname{Tot}}
\def\Diag{\operatorname{Diag}}

\def\set#1#2{\left\{\left. #1 \right| #2\right\}}

\newcommand{\N}{\mathbb N}

\numberwithin{equation}{subsection}
\newtheorem{lemma}[equation]{Lemma}
\newtheorem{thm}[equation]{Theorem}
\newtheorem{prop}[equation]{Proposition}
\newtheorem{defn}[equation]{Definition}
\newtheorem{remark}[equation]{Remark}

\title[Bialgebra cohomology]
{Bialgebra cohomology, pointed Hopf algebras, and deformations}
\author{Mitja Mastnak}
\thanks{The first author was partially supported by NSERC}
\address{Department of Pure Mathematics, University of Waterloo, 202
University
Ave West, Waterloo, ON N2L3G1, Canada}
\email{mmastnak@math.uwaterloo.ca}

\author{Sarah Witherspoon}
\thanks{The second author was partially supported by
the Alexander von Humboldt Foundation, NSF grant DMS-0443476,
and NSA grant H98230-07-1-0038.}
\address{Department of Mathematics, Texas A\&M University,
College Station, TX 77845, USA}
\email{sjw@math.tamu.edu}

\date{May 6, 2008}

\begin{document}

\maketitle

\begin{abstract}We give explicit formulas for maps in
a long exact sequence connecting bialgebra
cohomology to Hochschild cohomology.
We give a sufficient condition for the connecting homomorphism
to be surjective.
We apply these results to compute all
bialgebra two-cocycles of certain Radford biproducts
(bosonizations).
These two-cocycles are precisely those associated to
the finite dimensional pointed Hopf algebras in the recent
classification of Andruskiewitsch and Schneider, in an interpretation
of these Hopf algebras as graded bialgebra deformations of
Radford biproducts.
\end{abstract}

\section{Introduction}

Gerstenhaber and Schack \cite{GS2} found a long exact sequence connecting
bialgebra cohomology to Hochschild cohomology.
We reinterpret this sequence in the case of a finite dimensional Hopf
algebra:
Using results of Schauenburg \cite{Sc} and Taillefer \cite{T}, bialgebra
cohomology may be expressed as Hochschild cohomology of the Drinfeld
double;
we start by proving this directly in Section \ref{Dd}.
Other terms in the long exact
sequence (\ref{les}) involve Hochschild cohomology, with
trivial coefficients, of the Hopf algebra and of its dual.
This version of the sequence is particularly useful for computation.
We give an explicit description of the connecting homomorphism
in formulas (\ref{alphaf}) and (\ref{Gb}), and of the other two maps in
the long exact sequence  in Proposition \ref{restriction}
and formula (\ref{iota-bar}).
Theorem \ref{thm-main} gives a sufficient condition for the 
connecting homomorphism to be surjective in degree two.

As an application, we compute in Theorem \ref{Fbasis}
the (truncated) bialgebra cohomology, in
degree two, of the finite dimensional graded pointed Hopf algebras
arising in the recent classification of Andruskiewitsch and Schneider
\cite{AS8}.
They classified all finite dimensional pointed Hopf algebras having
abelian groups of grouplike elements, under a mild condition on
the group orders.
These include Lusztig's small quantum groups.
In general they are filtered and are
deformations (liftings) of their associated graded Hopf algebras,
an observation of Du, Chen, and Ye \cite{DCY} inspired by the
graded algebraic deformation theory of Braverman and Gaitsgory \cite{BG}.
These graded Hopf algebras are Radford biproducts, their deformations
governed by bialgebra cohomology in degree two.
In this setting, our Theorem \ref{thm-main} implies that 
the connecting homomorphism
in the long exact sequence (\ref{les}) is surjective (Theorem
\ref{conn-surj}).
We compute Hochschild cohomology with
trivial coefficients (Theorem \ref{H2B})
and apply the connecting homomorphism to give the degree two bialgebra
cohomology in Theorem \ref{Fbasis}.
This computation is analogous to that of
Gr\"unenfelder and the first author \cite{GM} of cohomology associated
to an abelian Singer pair of Hopf algebras.
We describe all homogeneous bialgebra two-cocycles of negative degree. It
follows from the classification of
Andruskiewitsch and Schneider that they all lift to deformations, providing
explicit examples for the Du-Chen-Ye theory.
We give a further
set of examples, the rank one pointed Hopf algebras of Krop and Radford
\cite{KR}, at the end of Section \ref{phad}.

Our computation of cohomology gives insight into the possible
deformations (liftings) of a Radford biproduct, providing
a different way to see {\em why} the liftings of Andruskiewitsch and
Schneider
must look the way they do.
In particular, we recover conditions for existence of certain
relations in the Hopf algebra from those for existence of
corresponding two-cocycles in Theorem \ref{H2B} and equations
(\ref{alphafij})
and (\ref{alphafii}).
Our computational techniques may be
useful in the search for pointed Hopf algebras that are left out of the
classification, that is those having small
prime divisors of the group order, complementing
work of Heckenberger \cite{H} on this problem.
These techniques should also be helpful in understanding
infinite dimensional pointed Hopf algebras.
We will address these problems in future papers.

We thank A.\ Masuoka, H.-J.\ Schneider and
P.\ Schauenburg for many helpful conversations; it was H.-J.\ Schneider
who first suggested this project to us.
We thank Ludwig-Maximilians-Universit\"{a}t M\"{u}nchen
for its hospitality during the preparation of this paper.

\section{Definitions and Preliminary Results}

All vector spaces (algebras, coalgebras, bialgebras)
will be over a ground field $k$.
In the classification of Andruskiewitsch and Schneider, $k$ must be
algebraically closed and of characteristic 0, however we do not
require this for the general theory.
If $A$ is an algebra and $C$ a coalgebra, then
$\Hom_k(C,A)$ denotes the convolution algebra of all linear maps from $C$
to $A$.
The unit and the multiplication on $A$ are denoted by $\eta\colon k\to A$
and
$m\colon A\ot A\to A$;
the counit and the comultiplication on $C$ are denoted by $\ep\colon C\to
k$ and
$\De\colon C\to C\ot C$.
We use Sweedler's notation for comultiplication: $\De(c)=c_1\ot c_2$,
($1\ot\De)\De(c)=c_1\ot c_2\ot c_3$, etc. If $f\colon U\ot V\to W$ is a
linear map,
then we often write
$f(u,v)$ instead of $f(u\ot v)$.
If $A$ is an augmented algebra, we denote the augmentation ideal
by $A^+=\ker\ep$. If $V$ is a vector space, we denote its $n$-fold
tensor power by $V^n$. To avoid confusion with comultiplication,
we write indices as superscripts, e.g. $v^1\ot\ldots\ot v^n\in V^n$.
If $A$ is an algebra, then we denote the $n$-ary multiplication by
$\widehat{(-)}$, i.e.\ if $\mathbf{a}=
a^1\ot\ldots\ot a^n\in A^n$, then $\widehat{\mathbf{a}}=a^1\cdots a^n$.
If $C$ is a coalgebra, then $\Delta^n$ denotes the $n$-ary
comultiplication, i.e.
$\Delta^n c=c_1\otimes\ldots\otimes c_n$.

\subsection{Bialgebra cohomology and deformations}\label{bcd}
We recall the definition of bialgebra cohomology and its truncated
version. For more details and greater generality we refer to \cite{GS}.
Let $B$ be a bialgebra. The left and right diagonal actions and coactions
of $B$ on
$B^n$ will be denoted by $\lambda_l, \lambda_r, \rho_l,\rho_r$,
respectively. More precisely, if
$a\in B$ and $\mathbf{b}=b^1\ot\ldots\ot b^n\in B^n$, then
\begin{eqnarray*}
\lambda_l(a\ot \b)&=&a_1b^1\ot\ldots\ot a_nb^n,\\
\lambda_r(\b\ot a)&=&b^1a_1\ot\ldots\ot b^na_n,\\
\rho_l(\b)&=&\widehat{\b}_1\ot\b_2=(b^1_1\ldots b^n_1)\ot
(b^1_2\ot\ldots\ot b^n_2),\\
\rho_r(\b)&=&\b_1\ot\widehat{\b}_2=(b^1_1\ot\ldots\ot b^n_1)\ot
(b^1_2\ldots b^n_2).
\end{eqnarray*}
The standard complex for computing bialgebra cohomology is the following
cosimplicial bicomplex
$\B^{p,q}$. The vertices are $\B^{p,q}=\Hom_k(B^p,B^q)$. The horizontal
faces
$$\na^h_i\colon \Hom_k(B^p, B^q)\to \Hom_k(B^{p+1}, B^q)$$ and
degeneracies
$$\sigma^h_i\colon \Hom_k(B^{p+1}, B^q)\to \Hom_k(B^p, B^q)$$
are those for computing Hochschild cohomology:
\begin{eqnarray*}
\na^h_{0} f&=&\lambda_l(1\ot f),\\
\na^h_{i} f&=& f(1\ot\ldots\ot m\ot\ldots \ot 1),\; 1\le i\le p,\\
\na^h_{p+1}f&=& \lambda_r(f\ot 1),\\
\sigma^h_{i} f&=& f(1\ot\ldots\ot\eta\ot\ldots\ot 1);
\end{eqnarray*}
the vertical faces
$$\na^c_j\colon \Hom_k(B^p, B^q)\to \Hom_k(B^p, B^{q+1})$$ and
degeneracies
$$\sigma^c_j\colon \Hom_k(B^p, B^{q+1})\to \Hom_k(B^p, B^{q})$$ are those
for computing
coalgebra (Cartier) cohomology:
\begin{eqnarray*}
\na^c_{0} f&=&(1\ot f)\rho_l,\\
\na^c_{j} f&=& (1\ot\ldots\ot \De\ot\ldots \ot 1)f,\; \ 1\le j\le q,\\
\na^c_{q+1}f&=& (f\ot 1)\rho_r,\\
\sigma^c_{i} f&=& (1\ot\ldots\ot\ep\ot\ldots\ot 1)f.
\end{eqnarray*}
The vertical and horizontal differentials are given by the usual
alternating
sums
$$\na^h=\sum (-1)^i\na^h_i,\; \ \ \ \na^c=\sum (-1)^j\na^c_j.$$
The {\em bialgebra cohomology} of $B$ is
then defined as
$$\coh_b^*(B)=\coh^*(\Tot\B).$$
where $$\Tot\B = \B^{0,0}\to \B^{1,0}\oplus\B^{0,1}\to\ldots \to
\bigoplus_{p+q=n} \B^{p,q}\stackrel{\partial^b}{\to}\ldots. $$
and $\partial^b$ is given by the sign trick (i.e.,
$\partial^b|_{\B^{p,q}}=\partial^h\oplus (-1)^1\partial^c\colon
\B^{p,q}\to \B^{p+1,q}\oplus\B^{p,q+1}$).
Here we abuse the notation by identifying a cosimplicial bicomplex with
its associated cochain bicomplex.
Let $\B_0$ denote the bicomplex obtained from $\B$ by replacing the edges
by zeroes, that is $\B_0^{p,0}=0=\B_0^{0,q}$ for all $p,q$.
The {\em truncated bialgebra cohomology} is
$$\widehat{\coh}_b^*(B)=\coh^{*+1}(\Tot\B_0).$$
For computations usually the normalized subcomplex $\B^+$ is used. The
normalized complex $\B^+$ is
obtained from the cochain complex $\B$ by replacing $\B^{p,q}=\Hom_k(B^p,
B^q)$ with the intersection
of degeneracies
$$(\B^+)^{p,q}=(\cap \Ker\sigma^h_i)\cap (\cap \Ker\sigma^c_j)\simeq
\Hom_k((B^+)^p,(B^+)^q).$$
Note that we can identify
$$
\widehat{\coh}_b^1(B)=\{f\colon B^+\to B^+|f(ab)=af(b)+f(a)b, \ \De
f(a)=a_1\ot f(a_2)+f(a_1)\ot a_2\}
$$
and
$$
\widehat{\coh}_b^2(B)=\widehat{{\rm Z}}_b^2(B)/\widehat{{\rm B}}_b^2(B),
$$
where
\begin{equation}\label{zhat}
\hspace{-.3in}\widehat{{\rm Z}}_b^2(B)\simeq \{(f,g)| \ 
   f\colon B^+\ot B^+\to B^+,\  g\colon B^+\to B^+\ot B^+,
\end{equation}

\vspace{-.26in}

\begin{eqnarray*}
\hspace{1.7in}&&af(b,c)+f(a,bc)=f(ab,c)+f(a,b)c, \\
&&c_1\ot g(c_2)+(1\ot\De)g(c)=(\De\ot 1)g(c)+g(c_1)\ot c_2, \\
&&f(a_1, b_1)\ot a_2b_2 -\De f(a,b)+a_1b_1\ot f(a_2, b_2) = \\
&&\hspace{1.4in}-(\De a)g(b)+g(ab)-g(a)(\De b) \}
\end{eqnarray*}
and
\begin{equation}\label{bhat}
\hspace{-2.1in}\widehat{{\rm B}}_b^2(B)\simeq \{(f,g)|\ 
   \exists h\colon B^+\to B^+,
\end{equation}

\vspace{-.26in}

\begin{eqnarray*}
\hspace{1.1in}&&f(a,b)=ah(b)-h(ab)+h(a)b \\
&&g(c)=-c_1\ot h(c_2)+\De h(c) -h(c_1)\ot c_2 \}.
\end{eqnarray*}

\vspace{.1in}

A {\em deformation} of the bialgebra $B$, over $k[t]$,
consists of a $k[t]$-bilinear
multiplication $m_t=m+tm_1+t^2m_2+\cdots$ and a comultiplication
$\Delta_t=\Delta+t\Delta_1+t^2\Delta_2+\cdots$ with
respect to which the $k[t]$-module
$B[t]:=B\ot_k k[t]$ is again  a bialgebra.
In this paper, we are interested only in those deformations for which
$\Delta_t=\Delta$, since the pointed Hopf algebras in the
Andruskiewitsch-Schneider classification have this property.
Given such a deformation of $B$, let $r$ be the smallest positive
integer for which $m_r\neq 0$ (if such an $r$ exists).
Then $(m_r,0)$ is a two-cocycle in $\widehat{{\rm Z}}^2_b(B)$.
Every nontrivial deformation is equivalent to one for which the
corresponding
$(m_r, 0)$ represents a nontrivial cohomology class \cite{GS}.
Conversely, given a positive integer $r$ and a two-cocycle
$(m',0)$ in $\widehat{{\rm Z}}^2_b(B)$,
$m+t^r m'$ is an associative multiplication on
$B[t]/(t^{r+1})$, making it into a bialgebra over $k[t]/(t^{r+1})$.
There may or may not exist $m_{r+1},m_{r+2},\ldots$
for which $m+t^rm'+t^{r+1}m_{r+1}+t^{r+2}m^{r+2}+\ldots$ makes
$B[t]$ into a bialgebra over $k[t]$.
(For more details on deformations of bialgebras, see \cite{GS}.)

\subsection{Graded bialgebra cohomology}\label{gbc}
Here we recall the definition of graded (truncated) bialgebra cohomology
\cite{DCY}.
If $B$ is a graded bialgebra, then $\B_{(l)}$ denotes the subcomplex of
$\B$ consisting of
homogeneous maps of degree $l$, more precisely
$$
\B^{p,q}_{(l)}=\Hom_k(B^p,B^q)_l=\{f\colon B^p\to B^q|f\mbox{ is
homogeneous of degree } l\}.
$$
Complexes $(\B_0)_{(l)}$, $\B^+_{(l)}$ and $(\B^+_0)_{(l)}$ are defined
analogously. The graded bialgebra
and truncated graded bialgebra cohomologies are then defined by:
\begin{eqnarray*}
\coh^*_b(B)_l&=&\coh^*(\Tot\B_{(l)})=\coh^*(\Tot\B^+_{(l)}),\\
\widehat{\coh}_b^*(B)_l&=&\coh^{*+1}(\Tot(\B_0)_{(l)})=\coh^{*+1}(\Tot(\B^+_0)_{(l)}).
\end{eqnarray*}
Note that if $B$ is finite dimensional, then
$$
\coh_b^*(B)=\bigoplus_l \coh_b^*(B)_l \ \mbox{ and } \
\widehat{\coh}_b^*(B)=\bigoplus_l \widehat{\coh}_b^*(B)_l.
$$

An {\em $r$-deformation} of $B$ is a bialgebra deformation
of $B$ over $k[t]/(t^{r+1})$ given by $(m^r_t,\Delta^r_t)$. Given a graded
bialgebra two-cocycle $(m',\Delta ')$ of $B$, in degree $-r$,
there exists an $r$-deformation, given by $(m+t^rm', \Delta + t^r \Delta
')$.
In this paper, we only consider $r$-deformations for which
$\Delta^r_t=\Delta$.

\begin{remark}\label{rdiff}\rm
(cf.\ \cite[Prop.\ 1.5(c)]{BG}, \cite{Gr}) Suppose that
$(B[t]/(t^r), \ m_t^{r-1}, \ \Delta_t^{r-1})$ is an
$(r-1)$-deformation, where
$$m_t^{r-1}=m+t m_1+\ldots+t^{r-1} m_{r-1} \ \mbox{ and } \
\Delta_t^{r-1}=\Delta+t\Delta_1+\ldots +t^{r-1}\Delta_{r-1}.$$
If
$$ D=(B[t]/(t^{r+1}), m_t^{r-1}+t^{r} m_{r},
\Delta_t^{r-1}+t^{r}\Delta_{r})$$
is an $r$-deformation, then
$$D'=(B[t]/(t^{r+1}), m_t^{r-1}+t^{r}
m'_{r},\Delta_t^{r-1}+t^{r}\Delta'_{r})$$
is an $r$-deformation if, and only if,
$(m'_{r}-m_{r},\Delta'_{r}-\Delta_{r})\in
\widehat{{\rm Z}}_b^2(B)_{-r}.$
Note also that if $(m'_{r}-m_{r},\Delta'_{r}-\Delta_{r})\in
\widehat{{\rm B}}_b^2(B)_{-r}$, then deformations $D$ and $D'$ are
isomorphic.
\end{remark}

\subsection{Coradically trivial and cotrivial cocycle pairs}
In this section we collect some preliminary results about cocycles
that will be needed in Section \ref{phad}.
The first lemma largely follows from the theory of relative
bialgebra cohomology \cite{GS}; however we did not
find a proof in the literature and so we include one for
completeness.
Let $B$ be a graded
bialgebra, and let $p\colon B\rightarrow B_0$ denote the canonical
projection.

\begin{lemma}\label{trivial-cotrivial}
 If
$\mathrm{char} k=0$ and $B_0$ is either a group algebra or the
dual of a group algebra, then every $(f,g)\in \Zb^2(B)$ is
cohomologous to a cocycle pair $(f',g')$ for which $f'|_{B_0\ot
B+B\ot B_0}=0$ and $(p \ot 1)g'=0=(1\ot p)g'$. If $f=0$ (resp.
$g=0$) then we can assume that also $f'=0$ (resp. $g'=0$).
\end{lemma}

We say that $f'$ (respectively $g'$) is {\em trivial} (respectively {\em
cotrivial}) on $B_0$ in case it satisfies the conclusion of the lemma.

\begin{proof} Let $t\in B_0$ be the left and right integral in $B_0$ such
that
$\ep(t)=1$. Note also that $t_1\ot S(t_2)=S(t_1)\ot t_2$.
Recall that for $a\in B_0$ we have $t_1\ot S(t_2)a=at_1\ot S(t_2)$. We now
proceed as follows.

\noindent\textbf{Step 1}: For each $f$, we will construct  $s=s_f\colon
B\to B$ such that
\begin{enumerate}
\item $\partial^h(s)|_{B_0\ot B}= f|_{B_0\ot B}$.
\item If $f|_{B\ot B_0}=0$, then $\partial^h(s)|_{B\ot B_0}= 0$.
\item If $g=0$, then $\partial^c(s)=0$.
\item If $(p\otimes 1)g=0$, then $(p\otimes 1)\partial^c(s)=0$.
\item If $(1\otimes p)g=0$, then $(p\otimes 1)\partial^c(s)=0$.
\end{enumerate}
Define $s=s_f\colon B\to B$ by $s(b)=t_1f(S(t_2),b)$. We claim that $s$
has the required properties:
\begin{enumerate}
\item
For $a\in B_0$ and $b\in B$ we compute
\begin{eqnarray*}
(\nah s)(a,b)&=& a s(b)-s(a b)+s(a)b\\
&=& at_1 f(S(t_2),b)-t_1f(S(t_2), a b)+t_1 f(S(t_2),a)b\\
&=& t_1 f(S(t_2)a, b) - t_1f(S(t_2), a b)+t_1 f(S(t_2),a)b\\
&=& t_1S(t_2)f(a,b) - t_1S(t_2)f(a,b) + t_1 f(S(t_2)a, b)\\
&& - t_1f(S(t_2), a b)+t_1 f(S(t_2),a)b\\
&=& f(a,b)-t_1\left(\nah f(S(t_2),a,b)\right) = f(a,b).
\end{eqnarray*}
\item
If $f$ is such that $f|_{B\ot B_0}=0$, then $(\nah s)|_{B\ot
B_0}=0$:
\begin{eqnarray*}
(\nah s)(b,a) &=& b s(a) - s(ba) + s(b)a \\
&=& bt_1 f(S(t_2),a) - t_1f(S(t_2),ba)+ t_1 f(S(t_2),b) a\\
&=& -t_1f(S(t_2),ba)+t_1f(S(t_2),b) a \\
&=& -t_1S(t_2) f(b,a)+t_1f(S(t_2)b,a) + t_1(\nah f(S(t_2),b,a))\\
&=& 0.
\end{eqnarray*}
\item
\begin{eqnarray*}
(\nac s)(b)&=& b_1\ot s(b_2)-\De s(b)+s(b_1)\ot b_2 \\
&=& b_1\ot t_1f(S(t_2), b_2) - \De t_1 f(S(t_2),b) +t_1
f(S(t_2),b_1)\ot b_2 \\
&=& (t_1\ot t_2)(\nac f)(S(t_3),b) \\
&=& -(t_1\ot t_2)(\nah g)(S(t_3),b).
\end{eqnarray*}
\item
\begin{eqnarray*}
(p\ot 1)(\nac s)(b) &=& -(p\ot 1)(t_1\ot t_2)(\nah g)(S(t_3),b) \\
&=&-(t_1\ot t_2)\big{[}(p(S(t_4)b_1)\ot S(t_3)b_2)(p\ot
1)g(S(t_5)b_3)\\
&&- (p\ot 1)g(S(t_3)b)\\
&& +(((p\ot 1)g(S(t_5)b_1))(p(S(t_4)b_2)\ot
S(t_3)b_3))\big{]}\\
&=& 0.
\end{eqnarray*}
\item A symmetric version of the computation above works.
\end{enumerate}

\noindent\textbf{Step 2.}  Define $s'=s'_f\colon B\to B$ a right-hand side
version of $s$ by
$s'(b)=f(b,t_1)S(t_2)$ and note that $s'$ has properties analogous to
those for $s$.
Now define $r=r_f\colon B\to B$ by $r_f = s'_f + s_{f-\partial^h s'_f}$
and observe that
$f_{B_0\ot B+ B\ot B_0}=\nah r_f$ and that $\nac r_f$ is $B_0$-cotrivial
(resp. equal to $0$) whenever
$g$ is such.

\noindent\textbf{Step 3.} We dualize Step 2. Note that $g^*\colon
B^*\otimes B^*= (B\otimes B)^*\to B^*$ is a Hochschild cocycle and
$r_{g^*}\colon B^*\to B^*$ (see Step 2) is such that
$\partial^h r_{g^*}|_{B_0^*\ot B^*+B^*\ot B_0^*}= g^*|_{B_0^*\ot
B^*+B^*\ot B_0^*}$ and $\partial^c r_{g^*}$ is $B_0^*$-cotrivial
(resp. equal to $0$) whenever $f^*$ is  $B_0^*$-cotrivial (resp.
equal to $0$). Now dualize again to obtain $u_g:=r^*_{g^*}\colon
B\to B$ and note that $g-\nac u_g$ is $B_0$-cotrivial and that
$\nah u_g$ is $B_0$-trivial (resp. equal to $0$) whenever $f$ is
such.

\noindent\textbf{Step 4.} Define $v=v_{f,g}\colon B\to B$ by
$v=u_g+s_{f-\nah u_g}$ and note that $(f',g'):=(f,g)-(\nah v,\nac
v)$ is a $B_0$-trivial, $B_0$-cotrivial cocycle pair.
\end{proof}

\begin{remark}{\em
The above proof shows that the conclusion of the Lemma
\ref{trivial-cotrivial} holds whenever $B_0$ is either a
commutative or cocommutative semisimple and cosemisimple Hopf
algebra (with no assumptions on $k$).}
\end{remark}

\begin{remark}{\em
If $B=R\# B_0$ as an algebra for some algebra $R$,
and $f\colon B\ot B\to B$ is a $B_0$-trivial
Hochschild cocycle, then $f$ is uniquely determined by its values
on $B^+\ot B^+$. More precisely, if $x,x'\in R$ and $h,h'\in B_0$,
then $f(xh, x'h')= f(x, {^ h{x'}})$.}
\end{remark}

\begin{defn}
$$\Zb^2(B)^+=\set{(f,g)\in\Zb^2(B)}{f\; \mbox{is}\; B_0\mbox{-trivial}, \
g\;\mbox{is}\; B_0\mbox{-cotrivial}}.$$
\end{defn}
If $f\colon B\ot B\to B$, and $r$ is a nonnegative integer, then
define $f_r\colon B\ot B\to B$
by $f_r|_{(B\ot B)_r}= f|_{(B\ot B)_r}$ and
$f_r|_{(B\ot B)_s}=0$ for $s\not=r$. If $g\colon B\to B\ot B$, then
we define $g_r$ analogously. Note that $f=\sum_{r\ge 0} f_r$ and
$g=\sum_{r\ge 0} g_r$. Define $f_{\le r}$ by  $f_{\le
r}=\sum_{0\le i\le r}f_i$ and then $f_{<r}$, $g_{\le r}$, $g_{<r}$
in similar fashion.

We will need the following lemma.

\begin{lemma}\label{ll1}
Let $r$ be a positive integer and let $f\colon B\ot B\to B$ be a
homogeneous Hochschild cocycle (with respect to the left and right regular
actions of $B$). If $f_{< r}=0$, then $f_r\colon B\ot B\to
B$ is an $\ep$-cocycle (i.e.\ a cocycle with respect to the trivial
action of $B$ on $B$).
\end{lemma}
\begin{proof}
We need to check that for homogeneous $x,y,z\in B$ with $deg(x)$,
$deg(y)$, $deg(z)>0$ we have $f_r(xy,z)=f_r(x,yz)$. Indeed, if
$deg(x)+deg(y)+deg(z)\not=r$, then both sides are equal to $0$. If
$deg(x)+deg(y)+deg(z)=r$, then note that $f(x,y)=0=f(y,z)$ and
hence $xf(y,z)-f(xy,z)+f(x,yz)-f(x,y)z=0$ gives the desired
conclusion.
\end{proof}

\begin{lemma} \label{gr0}
Assume $B$ is generated in degrees 0 and 1.
\begin{itemize}\item[(i)] If $(f,g)\in\Zb^2(B)$, $r>1$, $f_{\le r}=0$, and
$g_{<r}=0$,
then $g_r=0$.
\item[(ii)] If $(f,g)\in\Zb^2(B)_l$, $l<-1$, $r>0$, and $f_{\le r}=0$,
then $g_{\le r}=0$.
\item[(iii)] If $(0,g)\in \Zb^2(B)_l$ and $l<-1$, then $g=0$.
\item[(iv)] If $(f,g)\in\Zb^2(B)_l^+$, $l<0$, $r>0$, and $f_{\le r}=0$,
then $g_{\le r}=0$.
\item[(v)] If $(0,g)\in \Zb^2(B)_l^+$ and $l<0$, then $g=0$.
\end{itemize}
\end{lemma}
\begin{proof}
\begin{itemize}
\item[(i)] Note that $B_{r}$ is spanned by elements $xy$, where
$x\in B_1$ and $y\in B_{r-1}$. Now observe $ g_r(xy)=(\De
x)g_r(y)+g_r(x)(\De y)-(\nac f)(x,y) = 0+0-0=0. $ \item[(ii)] Note
the homogeneity of $g$ implies that $g_{\le -l-1}=0$. Hence if
$r\le -l+1$ then we are done. If $r>-l+1$, then use induction and
part (i). \item[(iii)] Follows from (ii). \item[(iv)] Note that
$B_0$-cotriviality of $g$ implies that $g_{\le -l+1}=0$. Hence if
$r\le -l+1$ then we are done. If $r>-l+1$, then use induction and
part (i). \item[(v)] Follows from (iv).
\end{itemize}
\end{proof}

Using notation similar to that for bialgebra cohomology, we define
the following in relation to Hochschild cohomology:
\begin{eqnarray*}
\HZ^2(B,B) &=& \set{f\colon B\ot B\to B}{\nah f=0},\\
\HB^2(B,B) &=& \set{\nah h}{h\colon B\to B}, \\
\HH^2(B,B) &=& \HZ^2(B,B)/\HB^2(B,B), \\
\HZ^2(B,B)_l &=& \set{f\in \HZ^2(B,B)}{f\mbox{ homogeneous of degree
}l},\\
\HB^2(B,B)_l &=& \set{f\in \HB^2(B,B)}{f\mbox{ homogeneous of
degree }l}\\ &=& \set{\nah h}{h \mbox{ homogeneous of degree }l}, \\
\HH^2(B,B)_l &=& \HZ^2(B,B)_l/\HB^2(B,B)_l.
\end{eqnarray*}

Define $\pi\colon \Zb^2(B)\to \HZ^2(B)$, by $\pi(f,g)=f$. Note
that $\pi$ maps subspaces $\Bb^2(B)$, $\Zb^2(B)_l$, and $\Bb^2(B)_l$
into $\HB^2(B)$, $\HZ^2(B)_l$, and $\HB^2(B)_l$ respectively. Hence
$\pi$ gives rise to a map from $\Hb^2(B)$ to $\HH^2(B,B)$ and a
map from $\Hb^2(B)_l$ to $\HH^2(B,B)_l$. We abuse notation by denoting
these maps by $\pi$ as well.

We have the following relation between truncated bialgebra cohomology
and Hochschild cohomology in degree two.

\begin{thm}
Assume $B$ is generated in degrees $0$ and $1$. If either $l<-1$
or $l=-1$ and $B_0$ is either a group algebra or a dual of a group
algebra, then $\pi\colon\Hb^2(B)_l\to\HH^2(B,B)_l$ is injective.
\end{thm}
\begin{proof}
Suppose $(f,g)\in \Zb^2(B)_l$ represents a cohomology class in
$\Hb^2(B)_l$ such that $f=\pi(f,g)\in\HB^2(B,B)$. Note that if
$l=-1$, then we can without loss of generality assume that
$(f,g)\in\Zb^2(B)_l^+$, hence $f=0$ and thus by Lemma
\ref{gr0}(v), also $g=0$. Now assume that $l<-1$. Since $f$ is a
Hochschild coboundary there is $s\colon B\to B$ such that $f=\nah
s$ and hence $(f,g)\sim (f,g)-(\nah s, \nac s) = (0,g-\nac s)$. By
Lemma \ref{gr0}(ii)  this means that $g-\nac s= 0$ and hence
$(f,g)\in\Bb^2(B)_l$.
\end{proof}

\subsection{Cocycles stable under a group action}
In this section we explain how to identify cocycles stable
under a group action with cocycles on a smash product with a
group algebra; this identification is useful in explicit computations
such as those in the last section.
Let $R$ be a $k$-algebra with an action of a finite group $\Gamma$
by automorphisms.
Let $R\# k\Gamma$ denote the corresponding {\em smash product algebra},
that is $R\# k\Gamma$ is a free left $R$-module with basis $\Gamma$
and algebra structure given by $(rg)(sh)= (r(g\cdot s))(gh)$ for
all $r,s\in R$, $g,h\in G$.
If the characteristic of $k$ does
not divide the order of $\Gamma$, then
\begin{equation}\label{hochiso}
   \coh^*_h(R\# k\Gamma, k)\simeq \coh^*_h(R,k)^{\Gamma}
\end{equation}
(see for example \cite[Cor.\ 3.4]{S}). 
Let $B=R\# k\Gamma$, which need not be a bialgebra in 
this subsection.
If $f\colon R\otimes R\to k$ is a $\Gamma$-stable
cocycle, then the corresponding cocycle $\bar{f}\colon B\otimes
B\to k$ is given by $\bar{f}(rg,r'g')=f(r, ^g r')$ for
all $r,r'\in R$, $g\in \Gamma$.
(We will use
the same notation $f$ in place of $\overline{f}$ for convenience).
This observation is a special case of the following general lemma (cf.\
\cite[Thm.~5.1]{CGW}):

\begin{lemma}\label{ft}
Let $f\in \Hom_k(R^n,k)\simeq
\Hom_{R^e}(R^{n+2},k)$ be a function representing an element of
$\HH^n(R,k)^{\Gamma}$ expressed in terms of the bar complex for
$R$. The corresponding function $\overline{f} \in \Hom_{k}(B^{ n},
k) \simeq \Hom_{B^e}(B^{ n+2}, k)$ expressed in terms of the bar
complex  for $B$ is given by
\begin{equation}
  \overline{f}(a_1 h_1\ot \cdots \ot a_n h_n) =
  f (a_1\ot {}^{h_1}a_2 \ot\cdots\ot
  {}^{h_1\cdots h_{n-1}} a_n)
\end{equation}
for all $a_1,\ldots ,a_n\in R$ and $h_1,\ldots,h_n\in \Gamma$.
\end{lemma}

\begin{proof} We sketch a proof for completeness; similar results
appear in \cite{CGW} and elsewhere for other choices of
coefficients. Let $\mathcal D = \oplus_{g\in \Gamma} R^e (g\otimes
g^{-1})$, a subalgebra of $B^e$. We claim that the bar resolution
for $B$ (as $B^e$-module) is induced from the $\mathcal
D$-projective resolution of $R$,
\begin{equation}\label{dres}
\cdots \stackrel{\delta_3}{\longrightarrow} {\mathcal D}_2
 \stackrel{\delta_2}{\longrightarrow} {\mathcal D}_1
 \stackrel{\delta_1}{\longrightarrow} {\mathcal D}_0
 \stackrel{m}{\longrightarrow} R \rightarrow 0,
\end{equation}
where ${\mathcal D}_0 = {\mathcal D}$ and
$${\mathcal D}_n =\Span_k \{ a_0h_0\ot \cdots \ot a_{n+1}h_{n+1} \mid
 a_i\in R , \ h_i\in \Gamma, \ h_0\cdots h_{n+1} = 1\}
$$
is a $\mathcal D$-submodule of $B^{\ot (n+2)}$. Indeed, a map
$\displaystyle{B^e\ot_{\mathcal D} {\mathcal D}_n \stackrel
{\sim}{\rightarrow} B^{\ot (n+2)}}$ is given by
$$
 (b_{-1}\ot b_{n+2})\ot (b_0\ot\cdots\ot b_{n+1}) \mapsto
   b_{-1}b_0\ot b_1\ot\cdots\ot b_n \ot b_{n+1}b_{n+2},
$$
and its inverse $\psi$ is
$$
  a_0h_0\ot a_1h_1\ot \cdots \ot a_{n+1}h_{n+1}\mapsto
  (1\ot h_0\cdots h_{n+1})\ot (a_0h_0\ot\cdots\ot a_{n+1} h_n^{-1}\cdots
   h_0^{-1}),$$
for $a_i\in R$ and $h_i\in \Gamma$.

There is a map $\phi$ from (\ref{dres}) to the bar complex for
$R$, as they are both $R^e$-projective resolutions of $R$,
$$
  a_0h_0\ot \cdots\ot a_{n+1}h_{n+1} \mapsto
   a_0\ot {}^{h_0}a_1\ot {}^{h_0h_1}a_2\ot\cdots \ot {}^{h_0\cdots
h_n}a_{n+1}.
$$
(See \cite[(5.2)]{CGW}.) Applying these maps $\psi, \phi$ of
complexes, together with the isomorphism
$\Hom_{B^e}(B^e\ot_{{\mathcal{D}}} {\mathcal{D}}_n,k)\simeq
\Hom_{\mathcal{D}}({\mathcal{D}}_n,k)$, we have
\begin{eqnarray*}
\overline{f}(a_1h_1\ot\cdots \ot a_nh_n)\!\! &=&\!\!
  \psi^*\phi^* f(1\ot a_1h_1\ot\cdots\ot a_nh_n\ot 1)\\
 \!\!&=&\!\! \phi^* f ((1\!\ot\! h_1\cdots h_n)\!\ot\! (1\!\ot\!
a_1h_1\!\ot\cdots\ot\! a_nh_n
   \!\ot\! h_n^{-1}\cdots h_1^{-1}))\\
 \!\! &=&\!\! \phi^* f (1\ot a_1h_1\ot\cdots\ot a_nh_n\ot h_n^{-1}\cdots
   h_1^{-1})\\
   \!\! &=&\!\! f (a_1\ot {}^{h_1}a_2 \ot  {}^{h_1h_2}a_3\ot
  \cdots\ot {}^{h_1\cdots h_{n-1}}a_n),
\end{eqnarray*}
since the image of $f$ is the trivial module $k$.
\end{proof}

\section{A Long Exact Sequence for Bialgebra Cohomology}
When we are dealing with a truncated double complex, a standard tool for
computing its cohomology is a long exact sequence. More precisely, if $\A$
is a
cochain bicomplex, $\A_0$ its truncated bicomplex and $\A_1$ its edge
bicomplex, then the
short exact sequence of cochain complexes
$$
0\to \Tot\A_0\to \Tot\A\to \Tot\A_1\to 0
$$
gives rise to a long exact sequence of cohomologies:
$$
\ldots\to \coh^*(\Tot\A_0)\to \coh^*(\Tot(\A))\to
\coh^*(\Tot(\A_1))\stackrel{\delta}{\to} \coh^{*+1}(\Tot(\A_0))\to\ldots,
$$
where the connecting homomorphism $\delta\colon \coh^*(\Tot(\A_1))\to
\coh^{*+1}(\Tot(\A_0))$ is induced by the
differential. In the context of bialgebra cohomology this was already used
in \cite{GS}.
Furthermore, if $\A$ is a cosimplicial bicomplex, then by the
Eilenberg-Zilber Theorem \cite{W} (see
\cite[Appendix]{GM} for the cosimplicial version) we have
$\coh^*(\Tot\A)\simeq \coh^*(\Diag\A)$. If $\A$ is associated to a pair of
(co)triples and a distributive law between them, then the cohomology of
$\Diag\A$ is the cohomology associated to the composed (co)triple. On the
other hand, if the bicomplex $\A$ arises from some mixed distributive law
then often one can, with some finiteness
assumptions, use some duality to ``unmix'' the distributive law. This
strategy worked
remarkably well when dealing with cohomology associated to an abelian
Singer pair of Hopf algebras
\cite{GM} and can also be applied to truncated bialgebra cohomology. For
the sake of simplicity we deal with
this aspect of theory on the level of (co)simplicial bicomplexes and do
not go into such generalities as
(co)triples and distributive laws between them.

\subsection{``Unmixed'' complex for computing bialgebra cohomology}
From now on assume that $B$ is a finite dimensional Hopf algebra.
Let $X=(B^\op)^*=(B^*)^\cop $. We will denote the
usual pairing $X\ot B\to k$ by $\p{\_}{\_}$, i.e.\
if $x\in X$ and $a\in B$, then $\p{x}{a} =x(a)$. Note
that $X$ and $B$ act on each other in the usual way (if $x\in X$ and
$a\in B$, then the actions
are denoted by $\l{a}{x}, x^a, \l{x}{a}, a^x$):
\begin{eqnarray*}
\p{\l{a}{x}}{b}=\p{x}{ba}; \ \ \l{a}{x}=\p{x_1}{a}x_2,\\
\p{x^a}{b}=\p{x}{a b}; \ \ x^a=\p{x_2}{a}x_1,\\
\p{y}{\l{x}{a}}=\p{y x}{a}; \ \  \l{x}{a}=\p{x}{a_2}a_1,\\
\p{y}{a^x}=\p{x y}{a}; \ \ a^x=\p{x}{a_1}{a_2}.
\end{eqnarray*}
Observe that the diagonal actions of $B$ on $X^n$ and of $X$ on $B^n$
are given by
\begin{eqnarray*}
\l{a}{\x}=\p{\widehat{\x}_1}{a}\x_2,\\
\x^a=\p{\widehat\x_2}{a}\x_1,\\
\l{x}{\a}=\p{x}{\widehat\a_2}\a_1,\\
\a^x=\p{x}{\widehat\a_1}{\a_2}.
\end{eqnarray*}
We use the natural isomorphism
$$\Homk(B^q,B^p)\simeq \Homk(X^p\ot B^q,k),$$
given by identifying linear maps $f\colon X^p\ot B^q\to k$ with linear
maps
$\bar{f}\colon B^q\to B^p$, by $f(\x\ot\b)=\p{\x}{\bar{f}(\b)}$, to obtain
a cosimplicial bicomplex
$$\mathbf{C}=\left(\Homk(X^p\ot B^q, k), (\nax)^*,
(\nab)^*\right)$$
from the complex $\B$ defined in Section \ref{bcd}.
The dual faces
\begin{eqnarray*}
\nab_i&=&(\nab_i)^{p,q}\colon X^{p}\ot B^{q+1}\to X^{p}\ot B^{q},\\
\nab_j&=&(\nab_j)^{p,q}\colon X^{p+1}\ot B^{q}\to X^{p}\ot B^{q}
\end{eqnarray*}
are
\begin{eqnarray*}
\nab_{0}(\x,\a)&=&(\x^{a^1},a^2\ot\cdots\ot a^{q+1})
=\p{\widehat\x_2}{a^1}(\x_1,a^2\ot\cdots\ot a^{q+1}),\\
\nab_{i}(\x,\a)&=&(\x,a^1\ot\cdots \ot a^ia^{i+1}\ot\cdots \ot a^{q+1}),\;
1\le i\le q,\\
\nab_{q+1}(\x,\a)&=&(\l{{a^{q+1}}}{\x},a^1\ot\cdots\ot a^{q})=
\p{\widehat\x_1}{a^{q+1}}(\x_2,a^1\ot\cdots\ot a^q),\\
\nax_{0}(\x,\a)&=&(x^2\ot\cdots \ot
x^{p+1},\a^{x^1})=\p{x^1}{\widehat\a_1}(x^2\ot\cdots \ot x^{p+1},\a_2)\\
\nax_{j}(\x,\a)&=&(x^1\ot\cdots\ot x^jx^{j+1}\ot\cdots\ot x^{p+1},\a),\;
1\le j\le p,\\
\nax_{p+1}(\x,\a)&=&(x^1\ot\cdots\ot
x^{p},\l{{x^{p+1}}}{\a})=\p{x^{p+1}}{\widehat\a_2}(x^1\ot\cdots\ot
x^{p},\a_1).
\end{eqnarray*}
The dual degeneracies
\begin{eqnarray*}
\sx_i&=&(\sx_i)^{p,q}\colon X^p\ot B^q\to X^{p+1}\ot B^q,\\
\sb_j&=&(\sb_j)^{p,q}\colon X^p\ot B^q\to X^{p}\ot B^{q+1}
\end{eqnarray*}
are given by
\begin{eqnarray*}
&&\sx_i(x^1\ot\cdots\ot x^p\ot\b)=x^1\ot\cdots\ot x^i\ot 1\ot
x^{i+1}\ot\cdots\ot x^p\ot\b,\\
&&\sb_j(\x\ot b^1\ot\cdots\ot b^q)=\x\ot b^1\ot\ldots\ot b^j\ot 1\ot
b^{j+1}\ot\cdots\ot b^{q},
\end{eqnarray*}
and the differentials
$$\nab\colon X^{p}\ot B^{q+1}\to X^{p}\ot B^{q},\;\;
 \ \ \nax\colon X^{p+1}\ot B^{q}\to X^{p}\ot B^{q}$$
are given by the usual alternating sums, i.e.
$$\nab=\sum (-1)^i\nab_{i},\;\; \ \
\nax=\sum (-1)^j\nax_{j}.$$

Note that by the cosimplicial version of the Eilenberg-Zilber Theorem we
have
$\coh^*(\Tot(\C))\simeq \coh^*(\Diag(\C))$.

\subsection{The diagonal complex and cohomology of the Drinfeld double}
\label{Dd}
Note that the differential
$$\na_d=(\na^d)^*\colon \Homk(X^n\ot B^n,k)\to \Homk(X^{n+1}\ot
B^{n+1},k)$$
in the diagonal complex $\Diag(\mathbf{C})$ is given by $(\na^d)^n
=\sum_{i=0}^{n+1}(-1)^k \na_i^d$, where $\na_i^d=\nax_i\nab_i$.

Recall that $D(B)=X\bowtie B$, the Drinfeld double of $B$, is $X\ot B$ as
coalgebra
and the multiplication is given by
$$(x\bowtie a)(y\bowtie b)=x
(\l{{a_1}}{{y^{S^{-1}(a_3)}})\bowtie
a_2 b= \p{y_1}{a_1}\p{y_3}{S^{-1}(a_3)}} xy_2\bowtie a_2 b.$$
The associated
flip $c\colon B\ot X\to X\ot B$, is given by
$$c(a,x)=\l{{a_1}}{x^{S^{-1}(a_3)}}\ot
a_2=\p{x_1}{a_1}\p{x_3}{S^{-1}(a_3)}x_2\ot a_2.$$
This map induces $c_{i,j}\colon B^i\ot X^j\to X^j\ot B^i$ and
$\tilde{c}_n\colon (X\bowtie B)^n\to X^n\ot B^n$ in the obvious way. Note
that $$c_{i,j}(\a\ot\x) = \p{\widehat\x_1}{\widehat\a_1}\p{\widehat
\x_3}{S^{-1}(\widehat\a_3)}\x_2\ot\a_2.$$
Define also a map $$\phi_n\colon X^n\ot B^n\to X^n\ot B^n$$ by
$$\phi_n(\x,\a)=\p{\widehat\x_1}{S^{-1}(\widehat\a_1)}\x_2\ot\a_2.$$

The following identities are due to the fact that in order to compute
$\tilde{c}_n$, we can apply $c$'s in arbitrary order.
\begin{eqnarray*}
\tilde{c}_{n+1}&=&(1\ot c_{1,n}\ot 1)(1\ot\tilde{c}_n)\\
\tilde{c}_{n+1}&=&(1\ot c_{n,1}\ot 1)(\tilde{c}_n\ot 1)\\
\tilde{c}_{i+j+1}&=&(1\ot c_{i,j}\ot 1)(1\ot c_{i,1}\ot c_{1,j}\ot
1)(\tilde{c}_i\ot 1\ot\tilde{c}_j)\\
\tilde{c}_{i+j+2}&=&(1\ot c_{i,j}\ot 1)(1\ot c_{i,2}\ot c_{2,j}\ot
1)(\tilde{c}_i\ot \tilde{c}_2\ot\tilde{c}_j)
\end{eqnarray*}
Recall that the standard complex for computing
$\coh^*_h(D(B),k)$, Hochschild cohomology of $D(B)$ with trivial
coefficients, is
given by
$$
\mathbf{D} : \ \ \ldots \to \Homk(D(B)^n,k)\stackrel{(\na^h)^*}{\to}
\Homk(D(B)^{n+1},
k)\to\ldots
$$
where $\na^h=(\na^h)^n=\sum_{i=0}^{n+1} (-1)^i\na^h_i$ and
\begin{eqnarray*}
\na^h_0 (u^1\ot\ldots\ot u^{n+1})&=&\ep(u^1)(u^2\ot\ldots \ot u^{n+1}),\\
\na^h_i (u^1\ot\ldots\ot u^{n+1})&=&u^1\ot\ldots \ot u^iu^{i+1}\ot\ldots
\ot
u^{n+1},\\
\na^h_{n+1}(u^1\ot\ldots\ot u^{n+1})&=&\ep(u^{n+1})(u^1\ot\ldots \ot
u^{n}).
\end{eqnarray*}

\begin{thm}
The map $\phi_n\tilde{c}_n\colon (X\bowtie B)^n\to X^n\ot B^n$ induces an
isomorphism of complexes and hence $\coh^*_b(B)\simeq
\coh^* (\Diag(\mathbf{C}))\simeq \coh^*_h(D(B),k)$.
\end{thm}

\begin{proof}
Note that $\psi_n=\phi_n\tilde{c}_n$ is a linear isomorphism (
$\phi_n^{-1}(\x,\a)=\p{\widehat\x_1}{\widehat\a_1}\x_2\ot \a_2$,
$c^{-1}(x,a)=\p{x_3}{a_3}\p{S^{-1}(x_1)}{a_1}a_2\ot x_2$). We will
show that for every $n$ and $0\le i\le n+1$ the diagram
$$
\begin{CD}
D(B)^{n+1} @> \partial^h_i >> D(B)^{n}\\
@V \psi_{n+1} VV              @V \psi_n VV\\
X^{n+1}\ot B^{n+1} @> \partial^d_i >> X^{n}\ot B^{n}
\end{CD}
$$
commutes.
We first deal with the case $i=0$. Note that
$$\psi_n\na^h_{0}=\phi_n\tilde{c}_n(\ep\ot
1)=(\ep\ot\phi_n)(1\ot\tilde{c}_n)$$
and that
$$\na^d_{0}\psi_{n+1}=\na^d_{0}\phi_{n+1}(1\ot c_{1,n}\ot
1)(1\ot\tilde{c}_n).$$
Hence it is sufficient to prove that
$$(\ep\ot\phi_n)=(\na^d_{0})^n\phi_{n+1}(1\ot c_{1,n}\ot 1).$$ This is
achieved
by the following computation
\begin{eqnarray*}
\lefteqn{[\na^d_{0}\phi_{n+1}(1\ot c_{1,n}\ot 1)]((x\bowtie a)\ot
\y\ot\b)}\\
&=&\na^d_{0}\phi_{n+1}\p{\widehat\y_1}{a_1}\p{\widehat\y_3}{S^{-1}(a_3)}(x\ot
\y_2, a_2\ot\b)\\
&=&\na^d_{0}\p{\widehat\y_1}{a_1}\p{\widehat\y_4}{S^{-1}(a_4)}\p{x_1\widehat\y_2}{S^{-1}(a_2\hb_1)}
(x_2\ot \y_3, a_3\ot\b_2)\\
&=&\p{\widehat\y_1}{a_1}\p{\widehat\y_5}{S^{-1}(a_5)}\p{x_1\widehat\y_2}{S^{-1}(a_2\hb_1)}
\p{x_3\widehat\y_4}{a_3}\p{x_2}{a_4\widehat\b_2}(\y_3, \b_3)\\
&=&\p{\widehat\y_1}{a_1}\p{\widehat\y_5}{S^{-1}(a_6)}\p{x_1}{S^{-1}(a_3\hb_2)}
\p{\widehat\y_2}{S^{-1}(a_2\hb_1)}\\
&&\cdot \p{x_3}{a_4}\p{\widehat\y_4}{a_5}\p{x_2}{\widehat\b_3}(\y_3,
\b_4)\\
&=&\p{x}{a_4\hb_3S^{-1}(\hb_2)S^{-1}(a_3)}\p{\widehat\y_1}{S^{-1}(\hb_1)S^{-1}(a_2)a_1}\p{\widehat\y_3}{S^{-1}(a_6)a_5}
(\y_2,\b_4)\\
&=&\ep(a)\ep(x)\p{\hy_1}{S^{-1}(\hb_1)}(\y_2,\b_2)\\
&=&(\ep\ot\phi_n)((x\bowtie a)\ot \y\ot\b).
\end{eqnarray*}
A similar computation applies to $i=n+1$. The remaining cases, where $1\le
i\le n$, are
settled by the diagram below (where each of the squares is easily seen to
commute).
{\small
$$
\begin{CD}
D(B)^{i-1}\!\ot\! D(B)^2\!\ot\! D(B)^{n-i}\!\!\! @>  \partial^h_i=1\ot
m\ot 1 >> \!\!\! D(B)^{i-1}\!\ot\! D(B)\!\ot\! D(B)^{n-i}\\
@V \tilde{c}_{i-1}\ot\tilde{c}_{2}\ot\tilde{c}_{n-i} VV
@V\tilde{c}_{i-1}\ot 1\ot\tilde{c}_{n-i}VV\\
X^{i-1}\!\ot\! B^{i-1}\!\ot\! X^2\!\ot\! B^2\!\ot\! X^{n-i}\!\ot\!
B^{n-i}\!\!\!\!\!\!\!\!
@> 1\ot 1\ot m\ot m\ot 1\ot 1 >>
\!\!\!\!\!\!\!\! X^{i-1}\!\ot\! B^{i-1}\!\ot\! X\!\ot\! B\!\ot\!
X^{n-i}\!\ot\! B^{n-i}\\
@V 1\ot {c}_{i-1,2}\ot {c}_{2,n-i}\ot 1 VV  @V 1\ot {c}_{i-1,1}\ot
{c}_{1,n-i}\ot 1 VV  \\
X^{i+1}\!\ot\! B^{i-1}\!\ot\! X^{n-i}\!\ot\! B^{n-i+2}\!\! @> (1\ot m)\ot
1\ot 1\ot (m\ot 1) >>
\!\! X^i\!\ot\! B^{i-1}\!\ot\! X^{n-i}\!\ot\! B^n \\
@V 1\ot c_{i-1,n-i}\ot 1 VV @V 1\ot c_{i-1,n-i}\ot 1 VV\\
X^{n+1}\!\ot\! B^{n+1}\!\! @> \partial^d_i=(1\ot m\ot 1)\ot (1\ot m\ot 1)
>> \!\! X^n\!\ot\! B^n \\
@V \phi_{n+1} VV @V \phi_{n} VV\\
X^{n+1}\!\ot\! B^{n+1}\!\! @> \phantom{=(1\ot m\ot
1)}\partial^d_i\phantom{\ot (1\ot m\ot 1)} >> \!\! X^n\!\ot\! B^n
\end{CD}
$$}
\end{proof}
\begin{remark}{\em
The isomorphism $\coh_b^*(B)\simeq \coh_h^*(D(B),k)$ can also be deduced
from a result of  Taillefer \cite{T},
combined with the fact due to Schauenburg \cite{Sc}
 that the category of Yetter-Drinfeld modules is equivalent to the
category
of Hopf bimodules.
See the remark following Proposition 4.6 in \cite{T}.
}\end{remark}

\subsection{Long exact sequence}
Let $\C_0$ denote the bicomplex obtained from $\C$ by replacing the edges
by zeroes and
let $\C_1$ denote the edge subcomplex of $\C$. Then we have a short exact
sequence
of bicomplexes
$$0\to \C_0 \to \C\to \C_1\to 0,$$
hence a short exact sequence of their
total complexes
$$0\to \Tot(\C_0)\stackrel{\iota}{\to} \Tot(\C)\stackrel{\pi}{\to}
\Tot(\C_1)\to 0,$$
which then gives rise to a long exact sequence of cohomologies
($i\ge 1$)
$$
\ldots\stackrel{\coh(\iota)}{\to}
\coh^i(\Tot(\C))\stackrel{\coh(\pi)}{\to}
\coh^i(\Tot(\C_1))\stackrel{\delta}{\to} \coh^{i+1}(\Tot(\C_0))\to\ldots
$$
Now use isomorphisms
\begin{eqnarray*}
&&\coh^i(\Tot(\C_1))\simeq \coh^i_h(B,k)\oplus \coh^i_h(X,k),\\
&&\coh^i(\Tot(\C))\simeq \coh^i(\Diag(\C))\simeq \coh^i_h(X\bowtie
B,k)\mbox{ and }\\
&&\coh^i(\Tot(\C_0))\simeq {\widehat{\coh}^{i-1}_b}(B),
\end{eqnarray*}
to get a long exact sequence (cf.\ \cite[\S8]{GS2})
\begin{equation}\label{les}
\ldots\stackrel{\bar\iota}{\to} \coh^i_h(D(B),k)
\stackrel{\bar{\pi}}{\to} \coh^i_h(X,k)\oplus \coh^i_h(B,k)
\stackrel{\delta}{\to} \widehat{\coh}_b^{i}(B)\to\ldots
\end{equation}

\subsection{Morphisms in the sequence}
Note that the morphism
$$\coh^i_h(B,k)\oplus \coh^i_h(X,k)\stackrel{\delta}{\to}
\widehat{\coh}_h^i(B)$$
corresponds to the connecting homomorphism in the long exact sequence and
is therefore given
by the differential, i.e.\ if $f\colon B^i\to k$ and $g\colon X^i\to k$
are cocycles, then $\de(f,g)=(\nax f,(-1)^i\nab g)$. More precisely
$$F:=\nax f\in
\Hom_k(B^i, B)\subseteq \bigoplus_{m+n=i+1}\Hom_k(B^m,B^n),$$
is given by
\begin{equation}\label{alphaf}
  F(\b)=f(\b_1)\widehat{\b_2}-f(\b_2)\widehat{\b_1}.
\end{equation}
If we identify $g$ with an element of $B^i$ ($g\in (X^i)^*\simeq
(B^i)^{**}\simeq B^i$),
then
$$G:=(-1)^i\nab g\in\Hom_k(B,B^i)\subseteq
\bigoplus_{m+n=i+1}\Hom_k(B^m,B^n)$$
is given by
\begin{equation}\label{Gb}
G(b)=(-1)^i\left( (\De^{i}b)g- g(\De^{i}b)\right).
\end{equation}
Recall that $\De^i b=b_1\ot\ldots\ot b_i$.

Using the cosimplicial Alexander-Whitney map, we can also show that the
map
$\overline{\pi}$ in the sequence (\ref{les}) above
is the double restriction:
\begin{prop}\label{restriction}
The map $$\coh^i_h(D(B),k)\stackrel{\bar\pi}{\to} \coh^i_h(B,k)\oplus
\coh^i_h(X,k)$$
is the restriction map in each component.
\end{prop}
\begin{proof}
We will establish the result by showing that the following diagram
commutes.
$$
\begin{CD}
\Tot^n(N\C) @> \Phi >> \Diag^n(N\C) @> (\phi\tilde{c})^* >> (D(B)^n)^*\\
@|  @. @V\res_2 VV\\
\Tot^n(N\C) @>\pi >> \Tot^n(N\C_1) @>\subseteq >> (X^n)^*\oplus (B^n)^*
\end{CD}
$$
Here $N\C$ denotes the normalized subcomplex of $\C$ (a map $f\colon
X^p\ot
B^q\to k$ is in $N\C$
if $f(x^1\ot\ldots\ot x^p,b^1\ot\ldots\ot b^q)=0$ whenever one of $x^i$ or
$b^j$ is a scalar) and
$\Phi$ denotes the Alexander-Whitney map (if
$f\in (X^p\ot B^q)^*\subseteq \oplus_{i+j=n} (X^i\ot B^j)^*$, then
$\Phi(f)\in (X^n\ot B^n)^*$
is given by $\Phi(f)=f \nax_{p+1}\ldots\nax_{n}\nab_0\ldots\nab_0$).
Note that
$\Phi(f)|_{X^n}=f(1_{X^p}\ot\ep_{X^{n-p}}\ot\eta_{B^q})$ and that
$\Phi(f)|_{B^n}=f(\eta_{X^p}\ot\ep_{B^{n-q}}\ot 1_{B^q})$. Hence, if $f$
is normal, then
$$\Phi(f)|_{X^n}=\begin{cases} f; &p=n\\
0; &p<n, \end{cases} \ \ \mbox{ and } \ \
\Phi(f)|_{B^n}=\begin{cases} f; &q=n\\
0; &q<n. \end{cases}$$
Also note that $\phi\tilde{c}|_{X^n}=1_{X^n}\ot \eta_{B^n}$ and
$\phi\tilde{c}|_{B^n}=\eta_{X^n}\ot 1_{B^n}$. Thus, if
$\mathbf{f}=(f_0,\ldots, f_n)\in
\bigoplus (X^i\ot B^{n-i})^*$ is a normal cocycle, then
$\res_{X^n}(\phi\tilde{c})^*\Phi\mathbf{f}(\mathbf{x})=
f_0(\mathbf{x})$ and
$\res_{B^n}(\phi\tilde{c})^*\Phi(\mathbf{f})(\mathbf{b})=f_n(\mathbf{b})$
and hence
$\res_2(\phi\tilde{c})^*\Phi\mathbf{f}=(f_0,f_n)=\pi(\mathbf{f})$.
\end{proof}

The map
$
{\widehat \coh}^n_b(B)\stackrel{\bar\iota}{\to}
\coh^{n+1}(D(B),k)
$
is given by the composite
\begin{eqnarray*}
{\widehat \coh}^n_b(B)&\stackrel{\simeq}{\to}&
\coh^{n+1}(\Tot(\C_1))\stackrel{\iota}{\to} \coh^{n+1}(\Tot(\C))\\
&\stackrel{\Phi}{\to}&
\coh^{n+1}(\Diag(\C))\stackrel{(\phi\tilde{c})^*}{\to} \coh^{n+1}(D(B),k).
\end{eqnarray*}
More precisely, if $\bar{f}\colon X^i\ot B^{n+1-i}\rightarrow k$
corresponds to $f\colon B^{n+1-i}\to B^i$, then
\begin{equation}\label{iota-bar}
\bar{\iota}f=\bar{f}\nax_{i+1}\ldots\nax_{n+1}\nab_0\ldots\nab_0\phi\tilde{c}.
\end{equation}

\subsection{Graded version}
Now assume that $B$ is a finite dimensional graded Hopf algebra. Note that
$X$ inherits the grading from
$B$ and is nonpositively graded, and $D(B)$ is graded by both positive and
negative integers.
Note that morphisms in the long exact sequence preserve degrees of
homogeneous maps and hence
for every integer $l$ we get a long exact sequence:
$$
\ldots\stackrel{\bar\iota}{\to} \coh^i_h(D(B),k)_l
\stackrel{\bar{\pi}}{\to} \coh^i_h(X,k)_l\oplus \coh^i_h(B,k)_l
\stackrel{\delta}{\to} \widehat{\coh}_b^{i}(B)_l\to\ldots
$$
Also note that if $l$ is negative, then $\coh_h^i(X,k)_l=0$ (as $X$ is
nonpositively graded and thus all homogeneous maps from $X$ to $k$ are of
nonnegative
degree), and hence in this case the sequence is
$$
\ldots\stackrel{\bar\iota}{\to} \coh^i_h(D(B),k)_l
\stackrel{\bar{\pi}}{\to} \coh^i_h(B,k)_l
\stackrel{\delta}{\to} \widehat{\coh}_b^{i}(B)_l\to\ldots
$$

\section{A sufficient condition for surjectivity of the connecting
homomorphism}\label{sufficient}

In this section, we give a sufficient condition for surjectivity of
the connecting homomorphism $\delta$ in degree 2 of the long
exact sequence (\ref{les}). The surjectivity will allow us to
compute fully the bialgebra cohomology in degree 2 for some
general classes of examples in the last section.

\subsection{Second Hochschild cohomology of a graded Hopf algebra
with trivial coefficients}\label{sm}
If $U\stackrel{f}{\to} V\stackrel{g}{\to} W$ is a sequence of
vector space maps such $gf=0$, then
$$
\frac{\ker f^*}{\im g^*} \simeq \left(\frac{\ker g}{\im
f}\right)^*\simeq \left[\ker\left(\tilde{g}\colon \frac{V}{\im
f}\to W\right)\right]^*,
$$
where $\tilde{g}$ is the map induced by $g$.
We apply this observation to an augmented algebra $R$ with augmentation
ideal $R^+$ and the map 
$$
R^+\otimes R^+\otimes R^+\stackrel{m\ot 1-1\ot m}
{\relbar\joinrel\relbar\joinrel\relbar\joinrel\relbar\joinrel\longrightarrow}
R^+\otimes R^+\stackrel{m}{\to} R^+
$$
to compute Hochschild cohomology of $R$ with trivial coefficients:
$$
\HH^2(R,k)\simeq \frac{\ker(m\ot 1-1\ot m)^*}{\im(m^*)}\simeq
\left[\ker\left( \tilde{m}\colon R^+\otimes_{R^+} R^+ \to
R^+\right)\right]^*,
$$
where we abbreviate $R^+\otimes_{R^+} R^+ = \frac{R^+\otimes
R^+}{\im(m\ot 1-1\ot m)}$ and $\tilde{m}$ is the map induced by
multiplication $m\colon R^+\otimes R^+\to R^+$. Also abbreviate
\begin{equation}\label{M-defn}
M:= \ker\left( \tilde{m}\colon R^+\otimes_{R^+} R^+ \to
R^+\right).
\end{equation}
The isomorphism above can be described explicitly
as follows. Choose $\phi\colon (R^+)^2=\Span\{xy | x,y\in R^+\}\to
R^+\otimes R^+$ a splitting of $m$. If we are given a linear map
$g\colon M\to k$, then define a cocycle $\bar{g}\colon R^+\ot
R^+\to k$ by $f=g(id - \phi m)$. If $f\colon R^+\ot R^+\to k$ is a
cocycle, then $\tilde{f}\colon M \to k$ is simply the induced map.
It is easy to check that $\tilde{\bar{g}}=g$ and that
$\bar{\tilde{f}}=f-\partial^h (f\phi)\sim f$.

\subsection{Surjectivity of the connecting homomorphism}
If $B=\bigoplus_{n\geq 0} B_n$ is a graded Hopf algebra, and
$p\colon B\to B_0$ is the canonical projection then $B$ equipped
with $B_0 \stackrel{p}{\leftrightarrows} B$ is Hopf algebra with a
projection in the sense of \cite{RA} and hence $R=B^{co
B_0}=\{r\in B | (1\otimes p)\Delta r = r\otimes 1\}$ is a Hopf
algebra in the category of Yetter-Drinfeld modules over $B_0$. The
action of $B_0$ on $R$ is given by $^h r = h_1 r S(h_2)$ and
coaction by $r\mapsto (p\otimes 1)\Delta r$. Throughout this
section we assume that $B_0=k\Gamma$ is a group algebra and that
the action of the group $\Gamma$ on $R$ is diagonal. In this case $R$ is
$(\Gamma\times \hat{\Gamma}\times\N)$-graded, that is it decomposes as
$R=\bigoplus R_{g,\chi, n}$, where $R_{g,\chi, n}$ consists of
homogeneous elements $r\in R$ of degree $n$ such that the coaction
of $k\Gamma$ is given by $r\mapsto g\otimes r$ and the action of
$k\Gamma$ is given by $^h r = \chi(h) r$.  We abbreviate
$R_{g,l}=\bigoplus_{\chi\in\hat\Gamma} R_{g,\chi,l}$.

Observe that $(m\ot 1-1\ot m)\colon R^+\ot R^+\ot R^+\to R^+\ot
R^+$ preserves the $(\Gamma\times \hat\Gamma\times
\mathbb{N})$-grading and hence we can decompose $M$ (see 
(\ref{M-defn})) in the same
fashion:
\begin{equation}\label{Mgxl}
M=\bigoplus_{(g,\chi,l)\in \Gamma\times\widehat{\Gamma}\times
\mathrm{Z}_{\ge 2}} M_{g,\chi,l},
\end{equation}
where $M_{g,\chi,l}$ consists
of homogeneous elements $m\in M$ of degree $l$ for which
the action and coaction of $k\Gamma$ are given by $^h m = \chi(h)
m$ and $m\mapsto g\otimes m$. Also note that
$\displaystyle{\HH^2(B,k)=\HH^2(R,k)^\Gamma\simeq
\bigoplus_{(g,n)\in \Gamma\times\mathbb{Z}_{\ge 2}}
M_{g,\varepsilon, n}^*}$ and that if $V$ is a finite-dimensional
trivial $B$-bimodule, then
$$\HH^2(B,V) = \bigoplus_{(g,n)\in \Gamma\times\mathbb{Z}_{\ge 2}}
\Hom_k(M_{g,\varepsilon, n}, V) = \bigoplus_{(g,n)\in
\Gamma\times\mathbb{Z}_{\ge 2}} M_{g,\varepsilon, n}^*\ot V.$$

The following lemma will be crucial in establishing a sufficient
condition for surjectivity of the connecting homomorphism.

\begin{lemma}\label{lm1} Let $B$ be a graded Hopf algebra of the
form described at the beginning of this section and let $l<0$.
Assume also that whenever $M_{h,\varepsilon,j}\not=0$ for some
$h\in\Gamma$ and $j>-l$, then $B$ contains no nonzero
$(1,h)$-primitive elements in degree $j+l$. If $R$ is generated as
an algebra by $R_1$, then for any $(f,g)\in \Zb^2(B)^+_l$ the
following holds:
\begin{itemize}
\item[(i)] If $r>-l$ and $f_{< r}=0$, then $(f,g)$ is cohomologous
to $(f',g')\in \Zb^2(B)^+_l$, where $f'_{\le r}=0$. \item[(ii)] If
$f_{-l}=0$, then $(f,g)\in \Bb^2(B)_l$.
\end{itemize}
\end{lemma}

\begin{proof}
(i) Note that by Lemma \ref{gr0}(iv) we have $g_{< r}=0$. By Lemma
\ref{ll1}, $f_r$ is an $\ep$-cocycle. If $u\in R^+\ot R^+$
represents an element in $M_{h,\chi,r}\not=0$, then, considering
(\ref{zhat}), we have
$$
0=g_r(0)=g_r(m(u))= h\ot f_r(u) - \De f_r(u) + f_r(u)\ot
1.
$$
Hence $f_r(u)=0$, since it is a $(1,h)$-primitive element of
degree $r+l$. Since $f_r(u)=0$ for all $u\in M$, we can conclude,
due to the discussion above, that $f_r$ is an
$\varepsilon$-coboundary. Thus we may let $s\colon B\to B$ be such
that $f_r(x,y)=s(xy)$ for $x,y\in B^+$. Note that $\nah s$ is
$B_0$-trivial, $\nac s$ is $B_0$-cotrivial and that
$(f',g')=(f,g)-(\nah s,\nac s)\in \Zb^2(B)_l^2$ is such that
$f'_{\le r}=0$.

(ii) Use induction and part (i) to show that $(f,g)\sim (0,g')$.
Then use Lemma \ref{gr0}(v).
\end{proof}

The following is one of the main theorems in our paper. The
results in the rest of the paper rely heavily on it.
Recall the notation defined in (\ref{Mgxl}).

\begin{thm}\label{thm-main}
Suppose that $l<0$ and $B$ is a finite dimensional graded Hopf
algebra such that
\begin{itemize}
\item $B$ is generated as an algebra by $B_0$ and $B_1$. \item
$B_0=k\Gamma$ and the action of $B_0$ on $R$ is diagonalizable,
i.e. $\Gamma$ acts on $R$ by characters. \item If
$M_{h,\varepsilon,j}\not=0$ for some $h\in\Gamma$ and $j>-l$, then
$B$ contains no nonzero $(1,h)$-primitive elements in degree
$j+l$.
\end{itemize}
Then the connecting homomorphism $\delta\colon \HH^2(B,k)_l\to
\Hb^2(B)_l$ is surjective.
\end{thm}
\begin{proof}
Let $(f,g)\in \Zb^2(B)^+_l$. Note that $f_{\leq - l -1}=0$ and
that $g_{\le -l+1}=0$. Define $\tilde f\colon B\ot B \to k$ by
$\tilde f(a,b)=-p_1f_{-l}(a,b)$, where $p_1\colon k\Gamma\to k$ is
given by $p_1(g)=\delta_{1,g}$. Note that $\tilde{f}$ is an
$\ep$-cocycle by Lemma \ref{ll1}. Now we prove that
$(f',g')=(f,g)-\de\tilde f=(f,g)-(\nac \tilde{f},0)\in
\Bb^2(B)_l$: This will follow from Lemma \ref{lm1}(ii) once we see
that $f'_{-l}=0$. Indeed, if $x\in R_{h_x,i}, \ y\in R_{h_y, -l-i}$,
then $f(x,y)\in B_0$ is $(h_x h_y,1)$-primitive by the same
argument as in the proof of Lemma \ref{lm1}(i). If $h_x h_y=1$,
then $f(x,y)$ is primitive and hence $0$. Otherwise $f(x,y)= a
(h_x h_y-1)$ for some $a\in k$. Note that $a=-p_1(a(h_x h_y-1))=
-p_1(f(x,y))= \tilde{f}(x,y)$ and that $\nac \tilde f(x,y)=
{\tilde f}(x_2,y_2)x_1y_1-{\tilde f}(x_1,y_1)x_2y_2=
\tilde{f}(x,y)(h_x h_y-1) = f(x,y)$. Thus $f'_{-l} = f_{-l} -\nac
\tilde{f} =0$.
\end{proof}

\begin{remark}{\em
If all $(1,h)$-primitive elements of $B$ in positive degree are
contained in $R_1$ (this happens whenever $B$ is coradically
graded and $B_0=k\Gamma$), then it is sufficient to demand that
there are no $(1,h)$-primitive elements in $R_1$, for all $h\in
\Gamma$ for which $M_{h,\varepsilon, l+1}\not=0$.}
\end{remark}

\section{Finite dimensional pointed Hopf algebras}\label{sec-fdpha}
We recall 
the Hopf algebras of Andruskiewitsch and Schneider \cite{AS8}, 
to which we will apply the results of the previous sections.

Let $\theta$ be a positive integer. Let $(a_{ij})_{1\leq i,j\leq
\theta}$ be a {\em Cartan matrix of finite type}, that is the
Dynkin diagram of $(a_{ij})$ is a disjoint union of copies of the
diagrams $A_{\bullet}, B_{\bullet}, C_{\bullet}, D_{\bullet}, E_6,
E_7, E_8, F_4, G_2$. In particular, $a_{ii}=2$ for $1\leq i\leq
\theta$, $a_{ij}$ is a nonpositive integer for $i\neq j$, and
$a_{ij}=0$ implies $a_{ji}=0$. Its Dynkin diagram is a graph with
vertices labelled $1,\ldots,\theta$. If $|a_{ij}|\geq |a_{ji}|$,
vertices $i$ and $j$ are connected by $|a_{ij}|$ lines, and these
lines are equipped with an arrow pointed toward $j$ if
$|a_{ij}|>1$.

Let $\Gamma$ be a finite abelian group. Let
$$
  {\mathcal D} = {\mathcal D}(\Gamma, (g_i)_{1\leq i\leq \theta},
   (\chi_i)_{1\leq i\leq \theta}, (a_{ij})_{1\leq i,j\leq \theta})
$$
be a {\em datum of finite Cartan type} associated to $\Gamma$ and
$(a_{ij})$; that is $g_i\in \Gamma$ and
$\chi_i\in\widehat{\Gamma}$ ($1\leq i\leq\theta$) such that
$\chi_i(g_i)\neq 1$ ($1\leq i\leq\theta$) and the Cartan condition
\begin{equation}\label{cartan}
   \chi_j(g_i)\chi_i(g_j)=\chi_i(g_i)^{a_{ij}}
\end{equation}
holds for $1\leq i,j\leq \theta$.

Let $\Phi$ denote the root system corresponding to $(a_{ij})$, and
fix a set of simple roots $\Pi$. If $\alpha_i, \alpha_j\in \Pi$,
write $i\sim j$ if the corresponding nodes in the Dynkin diagram
of $\Phi$ are in the same connected component. Choose scalars
$\lambda = (\lambda_{ij})_{1\leq i<j\leq \theta, \ i\not\sim j}$,
called {\em linking parameters}, such that
\begin{equation}\label{l-c}
   \lambda_{ij}=0 \ \ \mbox{ if } \ g_ig_j=1 \ \mbox{ or } \
\chi_i\chi_j\neq \varepsilon,
\end{equation}
where $\varepsilon$ is the trivial character defined by
$\varepsilon(g)=1$ ($g\in \Gamma$). Sometimes we use the notation
\begin{equation}\label{lambdaji}
   \lambda_{ji}:= -\chi_i(g_j)\lambda_{ij} \ \ \ (i<j).
\end{equation}
The (infinite dimensional) Hopf algebra $U({\mathcal D},\lambda)$
defined by Andruskiewitsch and Schneider \cite{AS8} is generated
as an algebra by $\Gamma$ and symbols $x_1,\ldots,x_{\theta}$,
subject to the following relations. Let $V$ be the vector space
with basis $x_1,\ldots,x_{\theta}$. The choice of characters
$\chi_i$ gives an action of $\Gamma$ by automorphisms on the
tensor algebra $T(V)$, in which $g(x_{i_1}\cdots
x_{i_s})=\chi_{i_1}(g)\cdots \chi_{i_s}(g)x_{i_1}\cdots x_{i_s}$
($g\in \Gamma$). We use this action to define the braided
commutators
$$
  \ad_c(x_i)(y) = [x_i,y]_c := x_i y - g_i(y) x_i,
$$
for all $y\in T(V)$. The map $c\colon T(V)\otimes T(V)\to
T(V)\otimes T(V)$, induced by $c(x_i\otimes y)=g_i(y)\otimes x_i$
is a braiding and $T(V)$ is a braided Hopf algebra in the
Yetter-Drinfeld category ${}^{\Gamma}_{\Gamma} {\mathcal {YD}}$.
(See \cite{AS8} for details, however we will not need to use the
theory of Yetter-Drinfeld categories.) There is a similar adjoint
action $\ad_c$ on any quotient of $T(V)$ by a homogeneous ideal.
The relations of $U({\mathcal D},\lambda)$  are those of $\Gamma$
and
\begin{eqnarray}
  gx_ig^{-1} &=& \chi_i(g) x_i  \ \ \ \ \ (g\in\Gamma, 1\leq i\leq
\theta),
\label{groupaction}\\
 (\ad_c(x_i))^{1-a_{ij}} (x_j) &=& 0 \ \ \ \ \  (i\neq j, \ i\sim j),
        \label{Serre}\\
 (\ad_c(x_i))(x_j)& =& \lambda_{ij} (1-g_ig_j)  \ \ \ \ \
    (i<j, \ i\not\sim j). \label{linking}
\end{eqnarray}
The coalgebra structure of $U({\mathcal{D}},\lambda)$ is defined
by
$$
 \Delta(g) = g\ot g , \ \ \ \Delta(x_i)=x_i\ot 1 + g_i\ot x_i,
$$
for all $g\in \Gamma$, $1\leq i\leq \theta$.

Let $W$ be the  Weyl group of the root system $\Phi$. Let
$w_0=s_{i_1}\cdots s_{i_p}$ be a reduced decomposition of the
longest element $w_0\in W$ as a product of simple reflections. Let
$$
  \beta_1=\alpha_{i_1}, \ \ \beta_2=s_{i_1}(\alpha_{i_2}),\ \cdots, \
 \beta_p=s_{i_1} s_{i_2}\cdots s_{i_{p-1}} (\alpha_{i_p}).
$$
Then $\beta_1,\ldots,\beta_p$ are precisely the positive roots
$\Phi^+$. Corresponding root vectors $x_{\beta_j}\in U({\mathcal
D},\lambda)$ are defined in the same way as for the traditional
quantum groups: In case $\mathcal D$ corresponds to the data for a
quantum group $U_q({\mathfrak{g}})$, let
$$
   x_{\beta_j}=T_{i_1}T_{i_2}\cdots T_{i_{j-1}}(x_{i_j}),
$$
where the $T_{i_j}$ are Lusztig's algebra automorphisms of
$U_q({\mathfrak{g}})$ \cite{L}. In particular, if $\beta_j$ is a
simple root $\alpha_l$, then $x_{\beta_j}=x_l$. The $x_{\beta_j}$
are in fact iterated braided commutators. In our more general
setting, as in \cite{AS8}, define the $x_{\beta_j}$ to be the
analogous iterated braided commutators.

The Hopf algebra $U({\mathcal D},\lambda)$ has the following
finite dimensional quotients. As in \cite{AS8} we make the
assumptions:
\begin{equation}\label{assumptions}
\begin{array}{l}
  \mbox{{\em the order of $\chi_i(g_i)$ is odd for all $i$,} }\\
  \mbox{{\em and is prime to 3 for all $i$ in a connected component of
type $G_2$.}}
\end{array}
\end{equation}
It follows that the order of $\chi_i(g_i)$ is constant in each
connected component $J$ of the Dynkin diagram \cite{AS8}; denote
this common order by $N_J$. It will also be convenient to denote
it by $N_{\alpha_i}$ or more generally by $N_{\beta_j}$ or $N_j$
for some positive root  $\beta_j$ in $J$. Let $\alpha\in\Phi^+$,
$\alpha=\sum_{i=1}^{\theta} n_i\alpha_i$, and let
$\Ht(\alpha)=\sum_{i=1}^{\theta}n_i$, $g_{\alpha}=\prod
g_i^{n_i}$, $\chi_{\alpha} =\prod \chi_i^{n_i}$. There is a unique
connected component $J_{\alpha}$ of the Dynkin diagram of $\Phi$
for which $n_i\neq 0$ implies $i\in J_{\alpha}$. We write
$J=J_{\alpha}$ when it is clear which $\alpha$ is intended. Choose
scalars $(\mu_{\alpha})_{\alpha\in\Phi^+}$, called {\em root
vector parameters}, such that
\begin{equation}\label{mu-alpha-c}
  \mu_{\alpha}=0 \ \mbox{ if } \ g_{\alpha}^{N_{\alpha}}=1 \
  \mbox{ or } \ \chi_{\alpha}^{N_{\alpha}}\neq \varepsilon.
\end{equation}
If $a=(a_1,\ldots, a_p)\in\N ^p - \{0\}$, define
$$
  \underline{a}:= a_1 \beta_1+\cdots +a_p\beta_p.
$$
In particular, letting $e_l:=(\delta_{kl})_{1\leq k\leq
p}\in\N^p-\{0\}$, we have $\underline{e_l}=\beta_l$.

The finite dimensional Hopf algebra $u({\mathcal D},\lambda,\mu)$
is the quotient of $U({\mathcal D},\lambda)$ by the ideal
generated by all
\begin{equation}\label{rootvector}
   x_{\alpha}^{N_{\alpha}} - u_{\alpha}(\mu) \ \ \ \ \ \ (\alpha\in
\Phi^+)
\end{equation}
where $u_{\alpha}(\mu)\in k\Gamma$ is defined inductively on
$\Phi^+$ as follows \cite[Defn.\ 2.14]{AS8}. If $\alpha$ is a
simple root, then $u_{\alpha}(\mu)
:=\mu_{\alpha}(1-g_{\alpha}^{N_{\alpha}})$. If $\alpha$ is not
simple, write $\alpha=\beta_l$ for some $l$, and then
\begin{equation}\label{ualpha}
   u_{\alpha}(\mu) := \mu_{\alpha}(1-g_{\alpha}^{N_{\alpha}}) +
   \sum_{\substack{b,c\in \N^p-\{0\} \\
\underline{b}+\underline{c}=\alpha}}
   t_{b,c}^{e_l} \mu_b u^c
\end{equation}
where
\begin{itemize}
\item[(i)] scalars $t^a_{b,c}$ are uniquely defined by
$$\Delta(x_{\beta_1}^{a_1N_1}\cdots x_{\beta_p}^{a_pN_p}) =
    x_{\beta_1}^{a_1N_1}\cdots x_{\beta_p}^{a_pN_p}\ot 1 +
   g_{\beta_1}^{a_1N_1}\cdots g_{\beta_p}^{a_pN_p}
   \ot  x_{\beta_1}^{a_1N_1}\cdots x_{\beta_p}^{a_pN_p}
$$

\vspace{-.2in}

$$\hspace{.5in} +\! \sum_{\substack{b,c\in\N^p-\{0\}  \\ \underline{b}+
  \underline{c} =\underline{a}}} \! t^a_{b,c} x_{\beta_1}^{b_1N_1}
  \cdots x_{\beta_p}^{b_pN_p} g_{\beta_1}^{c_1N_1}\cdots
g_{\beta_p}^{c_pN_p}
   \ot x_{\beta_1}^{c_1N_1}\cdots
     x_{\beta_p}^{c_pN_p} \ \ \ \mbox{\cite[Lemma 2.8]{AS8}}; $$
\item[(ii)] scalars $\mu_a$ and elements $u^a\in k\Gamma$ are
defined, via induction on $\Ht(\underline{a})$, by the
requirements that $\mu_{e_l}=\mu_{\beta_l}$ for $1\leq l\leq p$,
$\mu_a=0$ if $g_{\beta_1}^{a_1N_1}\cdots g_{\beta_p}^{a_pN_p}=1$,
and
$$
   u^a := \mu_a(1-g_{\beta_1}^{a_1N_1}\cdots g_{\beta_p}^{a_pN_p})
   + \sum_{\substack{b,c\in\N^p-\{0\} \\ \underline{b}+\underline{c}
  =\underline{a}}} t^a_{b,c} \mu_b u^c,
$$
where the remaining values of $\mu_a$ are determined by
$u^a=u^ru^s$ where $a=(a_1,\ldots,a_l,0,\ldots,0)$, $a_l\geq 1$,
$s=e_l$, and $a=r+s$ \ \ \cite[Theorem 2.13]{AS8}.
\end{itemize}

\noindent Andruskiewitsch and Schneider give the elements
$u_{\alpha}(\mu)$ explicitly in type $A_{\bullet}$ in
\cite[Theorem 6.1.8]{AS5}.

\begin{remark}\label{smallest}
{\em It follows from the induction \cite[Theorem 2.13]{AS8} that
if $\alpha$ is a positive root of {\em smallest} height for which
$\mu_{\alpha}\neq 0$, then $\mu_a=0$ for all $a\in\N ^p -\{0\}$
such that $\Ht(\underline{a}) < \Ht(\alpha)$. }\end{remark}

The following theorem is \cite[Classification Theorem 0.1]{AS8}.

\begin{thm}[Andruskiewitsch-Schneider]
Assume the field $k$ is algebraically closed and of characteristic
0. The Hopf algebras $u({\mathcal{D}}, \lambda,\mu)$ are finite
dimensional and pointed. If $H$ is a finite dimensional pointed
Hopf algebra having abelian group of grouplike elements with order
not divisible by primes less than 11, then $H\simeq
u({\mathcal{D}},\lambda,\mu)$ for some ${\mathcal D}$, $\lambda$,
$\mu$.
\end{thm}

We will need a lemma about central grouplike elements and skew
primitive elements.

\begin{lemma}\label{ll2}
Let $\alpha\in\Phi^+$ for which
$\chi_{\alpha}^{N_{\alpha}}=\varepsilon$. Then
\begin{itemize}\item[(i)] $g_{\alpha}^{N_{\alpha}}$ is in the
center of $U({\mathcal D},\lambda)$, and \item[(ii)] there are no
$(g_\alpha^{N_\alpha},1)$-skew primitives in $\bigoplus_{i\ge 1}
u({\mathcal D},\lambda,\mu)_i$.\end{itemize}
\end{lemma}

\begin{proof}
(i) It suffices to prove that $g_{\alpha}^{N_{\alpha}}$ commutes
with $x_j$ for each $j$. Note that  $g_{\alpha}^{N_{\alpha}} x_j =
\chi_j(g_{\alpha}^{N_{\alpha}}) x_j g_{\alpha}^{N_{\alpha}}$.
Write $g_{\alpha}=\prod  g_i^{n_i}$ where $\alpha =
\sum_{i=1}^{\theta} n_i\alpha_i$. By the Cartan condition
(\ref{cartan}) and the hypothesis
$\chi_{\alpha}^{N_{\alpha}}=\varepsilon$, we have
\begin{equation*}
  \chi_j(g_{\alpha}^{N_{\alpha}})  =
\prod_{i=1}^{\theta}\chi_j(g_i^{n_i})^{N_{\alpha}}
                = \prod_{i=1}^{\theta}(\chi_i(g_i)^{a_{ij}}
\chi_i^{-1}(g_j))^{n_iN_{\alpha}}
            = \chi_{\alpha}^{N_{\alpha}}(g_i^{a_{ij}}g_j^{-1}) = 1.
\end{equation*}

(ii) Each skew primitive in $u({\mathcal D},\lambda,\mu)$ is of
degree at most 1, and the only $(g,1)$-skew primitives in degree
1, for any $g\in\Gamma$, are in the span of
the $x_i$ \cite[(5.5) and Cor.\ 5.2]{AS8}. 
Now  $g_i\neq g_{\alpha}^{N_{\alpha}}$ for each $i$, since the latter element is
central by (i), while the former is not.
\end{proof}

We remark that the special case $u({\mathcal {D}},0,0)$ is a
graded bialgebra, the grading given by the coradical filtration.
In this case, $u({\mathcal {D}},0,0)\simeq {\mathcal{B}}(V)\#
k\Gamma$, the Radford biproduct (or bosonization) of the Nichols
algebra ${\mathcal{B}}(V)$ of the Yetter-Drinfeld module $V$ over
$k\Gamma$. For details, see \cite[Cor.\ 5.2]{AS8}.

We wish to understand $u({\mathcal D},\lambda,\mu)$ as a graded
bialgebra deformation of  $u({\mathcal D},0,0)$. We now describe
this graded case in more detail. Let $R={\mathcal{B}}(V)$ be the
subalgebra of $u({\mathcal D},0,0)$ generated by all $x_i$, and
$\widetilde{R}$ the subalgebra of $U({\mathcal D},0)$ generated by
all $x_i$, so that $R\simeq
\widetilde{R}/(x_{\alpha}^{N_{\alpha}}\mid \alpha\in\Phi^+)$. By
\cite[Thm.\ 2.6]{AS8}, $\widetilde{R}$ has PBW basis
\begin{equation}\label{PBW}
  x_{\beta_1}^{a_1}\cdots x_{\beta_p}^{a_p} \ \ \ \ \ \
  (a_1,\ldots,a_p\geq 0),
\end{equation}
and further,
\begin{equation}\label{central}
   [x_{\alpha},x_{\beta}^{N_{\beta}}]_c = 0
\end{equation}
for all $\alpha,\beta\in \Phi^+$. Thus $R$ has PBW basis
consisting of all elements in (\ref{PBW}) for which $0\leq a_i
<N_{i}$. Choose the section of the quotient map $\pi :
\widetilde{R}\rightarrow R$ for which the image of an element $r$
of $R$ is the unique element $\widetilde{r}$ that is a linear
combination of the PBW basis elements of $\widetilde{R}$ with
$a_i<N_{i}$ for all $i=1,\ldots,p$. This choice of section is used
in Section \ref{cna} below.

\section{Applications to some pointed Hopf algebras}\label{phad}

We will apply the cohomological results of the first part of the paper
to compute the degree 2 bialgebra cohomology
of the Radford biproduct $R\#
k\Gamma\simeq u({\mathcal{D}}, 0,0)$ defined in Section \ref{sec-fdpha}. 
We then use the result to understand deformations.

\subsection{Hochschild cohomology of $u({\mathcal{D}}, 0,0)$}\label{cna}
We first compute $\coh^2_h(R,k)$ and then apply the
isomorphism (\ref{hochiso}) to obtain $\coh^2_h(B,k)$ where $B=R\#
k\Gamma$. Hochschild one-cocycles on $R$ with coefficients in $k$
are simply derivations from $R$ to $k$, that is functions $f\colon
R\rightarrow k$ such that $f(rs)=\ep(r)f(s)+f(r)\ep(s)$ for all
$r,s\in R$. These may be identified with the linear functions from
$R^+/(R^+)^2$ to $k$, where $R^+=\ker\varepsilon$ is the
augmentation ideal. A basis for the space of such functions is $\{
f_i \mid 1\leq i\leq \theta\}$, where for each $i$,
$$
   f_i(x_j) = \delta_{ij} \ \ (1\leq j\leq \theta).
$$
All coboundaries in degree one are 0, and so
$\{f_i\mid 1\leq i\leq\theta\}$ may be identified with a basis of
$\HH^1(R,k)$.
We obtain some elements of $\HH^2(R,k)$ as cup products
of pairs of the $f_i$:
For $1\leq i<j\leq \theta$, define linear maps on pairs of PBW basis
elements
(\ref{PBW}), ${\bf x}^{\bf a} = x_{\beta_1}^{a_1}\cdots x_{\beta_p}^{a_p}$
and
${\bf x}^{\bf b}=x_{\beta_1}^{b_1}\cdots x_{\beta_p}^{b_p}$:
\begin{equation}\label{fij}
   f_{ji}({\bf x}^{\bf a},{\bf x}^{\bf b})=\left\{\begin{array}{rl}
    1, & \mbox{ if }{\bf x}^{\bf a}=x_{j} \mbox{ and }
            {\bf x}^{\bf b}=x_{i}\\
                   0, & \mbox{ otherwise.}\end{array}\right.
\end{equation}
Then $f_{ji}=f_j\smile f_i$.

Other Hochschild two-cocycles of $R$, with coefficients in $k$,
are indexed by the positive roots $\Phi^+$ and defined as follows:
Recall from the end of Section \ref{sec-fdpha} that $\widetilde{R}$ is an algebra
for which $R\simeq \widetilde{R}/(x_{\alpha}^{N_{\alpha}}\mid \alpha\in
\Phi^+)$.
Let $\widetilde{R}^+$ be the augmentation ideal of $\widetilde{R}$.
For each $\alpha\in\Phi^+$, define $\widetilde{f}_{\alpha}:
\widetilde{R}^+\ot \widetilde{R}^+ \rightarrow k$ by
$$
   \widetilde{f}_{\alpha}(r,
s)=\gamma_{(0,\ldots,0,N_{\alpha},0,\ldots,0)}
$$
where $N_{\alpha}$ is in the $i$th position if $\alpha=\beta_i$,
and $rs=\sum_{\bf a}\gamma_{\bf a}
{\bf x}^{\bf a}$ in $\widetilde{R}$.
By its definition, $\widetilde{f}_{\alpha}$ is associative on
$\widetilde{R}^+$,
so it may be extended (trivially) to a normalized Hochschild
two-cocycle on $\widetilde{R}$.
In fact $\widetilde{f}_{\alpha}$ is a coboundary on $\widetilde{R}$:
$\widetilde{f}_{\alpha}=\partial h_{\alpha}$ where
$h_{\alpha}(r)$ is the coefficient of $x_{\alpha}^{N_{\alpha}}$ in
$r\in\widetilde{R}$ written as a linear combination of PBW basis elements.
We next show that $\widetilde{f}_{\alpha}$ factors through the quotient
map $\pi :\widetilde{R}\rightarrow R$ to give a Hochschild two-cocycle
$f_{\alpha}$ on $R$, and that $f_{\alpha}$ is {\em not}
a coboundary on $R$. We must show that $\widetilde{f}_{\alpha}(r, s)=0$
whenever either $r$ or $s$ is in the kernel of the quotient map
$\pi:\widetilde{R}^+\rightarrow R^+$. It suffices to prove this for
PBW basis elements. Suppose ${\bf x}^{\bf a}\in\ker\pi$.
That is, $a_j\geq N_j$ for some $j$. Write ${\bf x}^{\bf a}=
\kappa x_{\beta_j}^{N_j} {\bf x}^{\bf b}$ where $\kappa$ is a nonzero
scalar and ${\bf b}$ may be 0; note this is possible by the
relation (\ref{central}). Then
$
  \widetilde{f}_{\alpha}({\bf x}^{\bf a},{\bf x}^{\bf c})=\kappa
  \widetilde{f}_{\alpha}(x^{N_j}_{\beta_j} {\bf x}^{\bf b},{\bf x}^{\bf
c}),
$
and this is the coefficient of $x_{\alpha}^{N_{\alpha}}$ in the product
$\kappa x_{\beta_j}^{N_j}{\bf x}^{\bf b}{\bf x}^{\bf c}$.
However, the coefficient of $x_{\alpha}^{N_{\alpha}}$ is 0:
If $\alpha=\beta_i$ and $j=i$, then since ${\bf x}^{\bf
c}\in\widetilde{R}^+$, this product
cannot have a nonzero coefficient for $x_{\alpha}^{N_{\alpha}}$.
If $j\neq i$, the same is true since $x_{\beta_j}^{N_j}$
is a factor of ${\bf x}^{\bf a} {\bf x}^{\bf c}$.
A similar argument applies to
$\widetilde{f}_{\alpha}({\bf x}^{\bf a},{\bf x}^{\bf c})$ if
${\bf x}^{\bf c}\in\ker\pi$.

Thus $\widetilde{f}_{\alpha}$ factors through $\pi\colon
\widetilde{R} \rightarrow R$, and we may define $f_{\alpha}\colon
R^+\ot R^+\rightarrow k$ by
\begin{equation}\label{fii}
   f_{\alpha}(r, s) = \widetilde{f}_{\alpha}(\widetilde{r},
    \widetilde{s}),
\end{equation}
where $\widetilde{r},\widetilde{s}$ are defined via the section of $\pi$
chosen at the end of Section \ref{sec-fdpha}.
We must verify that $f_{\alpha}$ is associative on $R^+$. Let $r,s,u\in
R^+$.
Since $\pi$ is an algebra homomorphism, we have
$\widetilde{r}\cdot \widetilde{s}=\widetilde{rs} + y$ and
$\widetilde{s}\cdot \widetilde{u}=\widetilde{su} + z$
for some elements $y,z\in\ker\pi$.
Since $\ker\pi\ot \widetilde{R}+\widetilde{R}\ot\ker\pi\subset
\ker\widetilde{f}_{\alpha}$, we have
\begin{eqnarray*}
   f_{\alpha}(rs, u) \ \ =\ \  \widetilde{f}_{\alpha}(\widetilde{rs},
   \widetilde{u})
   &=&  \widetilde{f}_{\alpha}(\widetilde{r}\cdot\widetilde{s}-y,
    \widetilde{u})\\
   &=&  \widetilde{f}_{\alpha}(\widetilde{r}\cdot\widetilde{s},
    \widetilde{u})\\
 &=&  \widetilde{f}_{\alpha}(\widetilde{r},\widetilde{s}\cdot
   \widetilde{u})
   \ \ =\ \   \widetilde{f}_{\alpha}(\widetilde{r},\widetilde{su})
   \ \ = \ \ f_{\alpha}(r,su).
\end{eqnarray*}

As we will see,
we  only need the functions $f_{ji}$ when $i\not\sim j$, that is
$i$ and $j$ are in different connected components of the Dynkin diagram
of $\Phi$.
Together with the $f_{\alpha}$, $\alpha\in\Phi^+$, these represent a
linearly independent subset of $\coh^2_h(R,k)$:

\begin{thm}\label{H2R}
The set $\{f_{\alpha}\mid \alpha\in\Phi^+\}\cup
\{f_{ji}\mid 1\leq i<j\leq\theta, \
i\not\sim j\}$ represents a linearly independent subset of
$\coh^2_h(R,k)$.
\end{thm}

\begin{proof}
Let $$f=\sum_{\alpha\in\Phi^+}c_{\alpha} f_{\alpha}
+ \sum_{\substack{1\leq i<j\leq\theta\\i\not\sim j}} c_{ji}f_{ji}$$
for scalars $c_{\alpha}, c_{ji}$.
Assume $f=\partial h$ for some $h:R\rightarrow k$.
Then for each $\alpha\in\Phi^+$,
$$
  c_{\alpha}=f(x_{\alpha},x_{\alpha}^{N_{\alpha}-1}) =\partial h
(x_{\alpha},
        x_{\alpha}^{N_{\alpha} -1}) = -h(x_{\alpha}^{N_{\alpha}}) =0
$$
as $x_{\alpha}, x_{\alpha}^{N_{\alpha} -1}\in R^+$ and
$x_{\alpha}^{N_{\alpha}}=0$
in $R$. For each pair $i,j$ ($1\leq i<j\leq \theta, \ i\not\sim j$),
$x_jx_i=\chi_i(g_j)x_ix_j$ since $i\not\sim j$, and so
\begin{eqnarray*}
  c_{ji} \ = \ f(x_j,x_i) & = & \partial h (x_j,x_i)\\
       & = & -h(x_jx_i)\\
   & = & -h(\chi_i(g_j)x_ix_j) \ = \ \chi_i(g_j)f(x_i,x_j) \ = \ 0.
\end{eqnarray*}
\end{proof}

Due to the isomorphism (\ref{hochiso}), we are primarily
interested in $\Gamma$-invariant Hochschild two-cocycles from
$R\ot R$ to $k$. The action of $\Gamma$ on the functions
$f_{\alpha},f_{ji}$ is diagonal, and so we determine those
$f_{\alpha}, f_{ji}$ that are themselves $\Gamma$-invariant.

\begin{thm}\label{H2B}
If $|\Gamma|$ is not divisible by primes less than $11$, then
$$
  \{f_{\alpha}, f_{ji}\}^\Gamma:=\{f_{\alpha}\mid \alpha\in \Phi^+, \ \chi_{\alpha}^{N_{\alpha}}
  =\varepsilon\}\cup\{f_{ji}\mid 1\leq i<j\leq\theta, \ i\not\sim j,
   \ \chi_{i}\chi_{j}=\varepsilon\}
$$
is a basis of $\coh^2_h(R\# k\Gamma,k)$.
In particular, if $l\leq -3$, then
$\{f_{\alpha} |  \chi_{\alpha}^{N_{\alpha}}=\varepsilon, \
  N_{\alpha}\Ht(\alpha) = -l\}$ is a basis for $\coh^2_h(R\# k\Gamma ,
k)_l$.
\end{thm}

\begin{remark}\label{conds}{\em
The condition $i\not\sim j$ is implied by the condition
$\chi_i\chi_j=\varepsilon$, and so is redundant:
A proof of this fact consists of a case-by-case
analysis using the Cartan condition (\ref{cartan}) for the
pairs $i,j$ and $j,i$.
We put the condition $i\not\sim j$ in the statement of the theorem
for clarity.
}\end{remark}

\begin{proof}
The action of $\Gamma$ on $\coh^2_h(R,k)$ comes from the dual
action of $\Gamma$ on $R\ot R$, that is, $(g\cdot f)(r,s)=f(g^{-1}\cdot r,
g^{-1}\cdot s)$. Therefore $g\cdot f_{\alpha}=\chi_{\alpha}^
{-N_{\alpha}}(g)f_{\alpha}$
and $g\cdot f_{ji}=\chi_{i}^{-1}(g)\chi_{j}^{-1}(g)
f_{ji}$.
Thus the subset of those functions from Theorem \ref{H2R} that are
$\Gamma$-invariant consists of the $f_{\alpha}$ for which
$\chi_{\alpha}^{N_{\alpha}}=\varepsilon$ and the $f_{ji}$
for which $\chi_i\chi_j=\varepsilon$.
By Remark \ref{conds},
this proves that the given set represents a linearly independent subset
of $\coh^2_h(R\# k\Gamma,k)$.

The fact that the given set spans $\coh^2_h(R\# k\Gamma,k)$
is a consequence of \cite[Lemma 5.4]{GK}:
That lemma states in our case that $\coh^2(R,k)$ (equivalently,
$\coh^2_h(R,k)$) has basis in one-to-one correspondence with that of
$I/(T^+(V)\cdot I + I\cdot T^+(V))$, where $R=T(V)/I$ ($I$ is the
ideal of relations).
This follows by looking at the minimal resolution of $k$.
As a consequence, $\coh^2_h(R\# k\Gamma,k)$, which is isomorphic
to $\coh^2_h(R,k)^{\Gamma}$, has basis in one-to-one correspondence
with that of the $\Gamma$-invariant subspace of $I/(T^+(V)\cdot I
+ I\cdot T^+(V))$.
The $\Gamma$-invariant root vector relations and linking relations
give rise to the elements $f_{\alpha}$ and $f_{ji}$ in the statement
of the theorem.
It remains to show that the Serre relations do not give rise to
$\Gamma$-invariant elements in cohomology.
If the Serre relation (\ref{Serre}) did give rise to a $\Gamma$-invariant
element in cohomology, then $\chi_i^{1-a_{ij}}\chi_j=\varepsilon$
by considering the $\Gamma$-action.
This is not possible:
The Cartan condition for the pairs $i,j$ and $j,i$, together with
this equation, implies $\chi_i(g_i)^{a_{ij}+a_{ji}-a_{ij}a_{ji}}=1$.
A case-by-case analysis shows that this implies $\chi_i(g_i)$ has
order 3, 5, or 7, contradicting our assumption that the order
of $\Gamma$ is not divisible by primes less than 11.

The last statement of the theorem is now immediate from the observation
that $f_{\alpha}$ is homogeneous of degree $-N_{\alpha}\Ht(\alpha)$, and
$f_{ji}$ is homogeneous of degree $-2$.
\end{proof}

\begin{remark}{\em
Masuoka  independently obtained a proof that the given set spans
$\coh^2_h(R\# k\Gamma,k)$,
using completely different methods and results from his preprint
\cite{Mas2}.}
\end{remark}

Note that the conditions in the theorem on the $\chi_{\alpha}$  and
$\chi_i$
are ``half'' of the conditions (\ref{l-c}) and (\ref{mu-alpha-c})
under which nontrivial linking
or root vector relations may occur.
The other half of those conditions, involving elements of $\Gamma$,
will appear after we apply the formula (\ref{alphaf}) to obtain
corresponding {\em bialgebra} two-cocycles. (The bialgebra two-cocycle
will be 0 when the
condition on the appropriate group element is not met.)

\subsection{Bialgebra two-cocycles}
Let $B=R\# k\Gamma$ as before. 
We wish to apply the connecting homomorphism
in the long exact sequence (\ref{les}) to elements of
$\coh^2_h(B,k)$ from Theorem \ref{H2B}, in order to obtain
bialgebra two-cocycles. First we prove that the connecting
homomorphism is surjective.

\begin{thm}\label{conn-surj} Assume the order of $\Gamma$ is not
divisible by $2$ or $3$. If $\HH^2(B,k)=\Span\{ f_{ji},
f_{\alpha}\}^\Gamma$ and $l<0$, then the connecting homomorphism
$\delta \colon \HH^2(B,k)_l\to\Hb^2(B)_l$ is surjective.
\end{thm}
\begin{proof}
Let $M$ be as in (\ref{M-defn}). Note that the fact that
$\HH^2(B,k)=\Span\{ f_{ji}, f_{\alpha}\}^\Gamma$ translates into
the fact that $\{x_j\otimes x_i - \chi_j(g_i) x_i\otimes x_j |
\chi_i\chi_j\}\bigcup\left\{x_\alpha\otimes x_\alpha^{N_\alpha-1}
| \chi_\alpha^{N_\alpha}=\varepsilon\right\}$ forms a basis for
$M^\Gamma$. Hence, if $M_{h,\varepsilon,r}\not=0$, then either
$h=g_i g_j$ with $\chi_i\chi_j = \varepsilon$, or
$h=g_\alpha^{N_\alpha}$ with $\chi_\alpha^{N_\alpha}=\varepsilon$.
If $h=g_\alpha^{N_\alpha}$, then there are no $(1,h)$-primitives
in of positive degree in $B$ by Lemma \ref{ll2}(ii). We now show that
there are no $(1,h)$-primitives of positive degree in $B$
(equivalently $R_1$) whenever $h=g_i g_j$. Suppose otherwise,
i.e., for some $k$ we have $g_k = g_i g_j$. Since $i\not\sim j$,
we also have that either $i\not\sim k$, or $j\not\sim k$. Without
loss of generality suppose the latter. Then
$\chi_k(g_k)=\chi_k(g_i g_j)\chi_i(g_k)\chi_j(g_k) =
(\chi_k(g_i)\chi_j(g_k))(\chi_k(g_j)\chi_j(g_k))=\chi_k(g_k)^{a_{ki}}$
and hence $\chi_k(g_k)^{1-a_{ki}}=1$. This is impossible, since
$1-a_{ki}\in \{1,2,3,4\}$ and $|\Gamma|$ is not divisible by $2$
or $3$.

Hence the conditions of Theorem \ref{thm-main} are satisfied.
\end{proof}

\begin{remark} {\em Assume that the order of $\Gamma$ is not divisible
by $2$ or $3$. Additionally assume that the order of $\Gamma$ is
not divisible by $5$ whenever the Dynkin diagram associated to
$\mathcal{D}$ contains a copy of $B_n$ with $n\ge 3$. A similar
case by case analysis to the one in the proof of Theorem \ref{H2B}
shows that there are no $(1,h)$-primitives whenever
$h=g_i^{1-a_{ij}} g_j$ with
$\chi_{i}^{1-a_{ij}}\chi_j=\varepsilon$. This observation can then
be used to show that the connecting homomorphism $\delta\colon
\HH^2(B,k)_l\to \Hb^2(B)_l$, where $B=u(\mathcal{D},0,0)$ and
$l<0$ is surjective.}
\end{remark}

Now let $f\colon B\ot B\rightarrow k$ be a Hochschild two-cocycle.
The formula (\ref{alphaf}) applied to $f$ yields
\begin{equation}\label{F2}
  F (a,b) = f(a_1,b_1)a_2b_2-f(a_2,b_2)a_1b_1,
\end{equation}
a bialgebra two-cocycle representing an element in
$\widehat{\coh}^2_b(B,B)$.
We apply this formula
to $f=f_{ji}, f_{\alpha}$, defined in (\ref{fij}) and (\ref{fii}),
to obtain explicit bialgebra two-cocycles $F=F_{ji}, F_{\alpha}$.
For our purposes, it will suffice to compute the value of each $F_{ji},
F_{\alpha}$ on a single well-chosen pair of elements in $R$.
In order to compute them on arbitrary pairs
of elements of $B$, one must use (\ref{F2}) and Lemma \ref{ft}.

\begin{lemma}
If $f\colon B\ot B\to k$ is a homogeneous $\varepsilon$-cocycle of
degree $l<0$, $\delta f = (F,0)$, and $x\in R_i$, $y\in R_j$, with
$i+j=-l$, are PBW-basis elements in components $g_x$, $g_y$, then
$$
F(x,y)=f(x,y)(1-g_xg_y).
$$
In particular
\begin{equation}\label{alphafij}
  F_{ji} (x_j,x_i) = 1-g_jg_i,
\end{equation}
and
\begin{equation}\label{alphafii}
  F_{\alpha}(x_{\alpha},x_{\alpha}^{N_{\alpha}-1}) =
1-g_{\alpha}^{N_{\alpha}}.
 \end{equation}
\end{lemma}
\begin{proof}
The hypotheses imply that
$\Delta x = x\ot 1+g_x\ot x + u$ and $\Delta y = y\ot 1 + g_y\ot y + v$,
where $u=\sum_r u'_r\ot u''_r\in \bigoplus_{p=1}^{i-1} B_p\ot B_{i-p}$
and
$v=\sum_s v'_s\ot v''_s\in \bigoplus_{q=1}^{j-1} B_q\ot B_{j-q}$. Then
$f(u'_r,v'_s)=0=
f(u''_r,v''_s)$ due to degree considerations and hence
$F(x,y)=f(x,y)1-f(x,y)g_xg_y$.
\end{proof}

Note that $F_{ji}(x_j,x_i)=0$ exactly when $g_ig_j=1$, and
$F_{\alpha}(x_{\alpha},x_{\alpha}^{N_{\alpha}-1})=0$ exactly when
$g_{\alpha}^{N_{\alpha}}=1$.  More generally, if $x,y$ are
PBW-basis elements of joint degree $2$ (resp. $\Ht(\alpha)N_\alpha$),
then $F_{ji}(x,y)\in k(1-g_jg_i)$ (resp. $F_\alpha(x,y)\in
k(1-g_\alpha^{N_\alpha})$).

Combined with the conditions on
$\chi_{\alpha}$
and $\chi_i$ in Theorem \ref{H2B}, we have recovered precisely
the conditions in (\ref{l-c}) and (\ref{mu-alpha-c}) under
which there exist nontrivial linking and root vector relations.
In Theorem \ref{thm-def} below, we make the connection
between these bialgebra two-cocycles and the pointed Hopf
algebras $u({\mathcal{D}}, \lambda,\mu)$.

Our calculations above, combined with
Theorems \ref{H2B} and \ref{conn-surj} now allow us
to determine completely $\Hb^2(B)_l$, $l<0$,
for coradically graded Hopf algebras $B=R\#\Gamma$ in the
Andruskiewitsch-Schneider program.
We have the following theorem.

\begin{thm}\label{Fbasis}
Let $B=u(\mathcal{D},0,0)$ and assume that $|\Gamma|$ is not
divisible by $2$ or $3$. If
$\HH^2(B,k)=\Span\{f_{ji},f_\alpha\}^\Gamma$ (see for example
Theorem \ref{H2B}), then the set
$$\set{(F_\alpha,0), (F_{j,i},0)}{\alpha\in\Phi^+, 1\leq i<j\leq
\theta, g_\alpha^{N_{\alpha}}\not=1, \chi_{\alpha}^{N_\alpha}=\ep,
i\not\sim j, g_i g_j \not=1, \chi_i\chi_j=\ep}$$ is a basis for
$\bigoplus_{l<0}\Hb^2(B)_l$.
\end{thm}
\begin{proof} $F_\alpha$ is a homogeneous cocycle of degree
$l:=-\Ht(\alpha)N_\alpha$.
If $g_\alpha^{N_\alpha}=1$, then $(F_\alpha)_{-l}=0$ and hence by Lemma
\ref{lm1} (ii)
$(F_\alpha,0)$ is a coboundary. We can show that if $g_i g_j=1$, then
$(F_{i,j},0)$ is a
coboundary in a similar fashion. Hence by Theorems \ref{H2B} and
\ref{conn-surj} and our above
calculations, the given set spans $\bigoplus_{l<0}\Hb^2(B)_l$.

Note that it is sufficient to show that for every $l<0$, the Hochschild
cohomology classes of
cocycles $F_\alpha, F_{i,j}$ of degree $l$ are linearly independent. This
is achieved as follows.

If $l<-2$, then $\left(\sum_{\Ht(\alpha)N_\alpha=-l}
\lambda_{\alpha}F_{\alpha}\right)(x_\alpha,x_\alpha^{N_\alpha-1})=\lambda_{\alpha}(1-g_\alpha^{N_\alpha})$,
but for
$s\colon B\to B$ homogeneous of degree $l$, we have
$\partial^h(s)(x_\alpha,x_\alpha^{N_\alpha-1})=x_\alpha
s(x_\alpha^{N_\alpha-1})-s(x_\alpha^{N_\alpha})+s(x_\alpha)x_\alpha^{N_\alpha-1}=0$.

If $l=-2$, then $\left(\sum_{i,j}
\lambda_{ji}F_{ji}\right)(x_j\otimes x_i- q_{i,j} x_i\otimes
x_j)=\lambda_{i,j}(1-g_ig_j)$ and if $s\colon B\to B$ is
homogeneous of degree $-2$, then $\partial^h(s)(x_j\otimes x_i -
\chi_j(g_i)x_i\otimes x_j)=0$.
\end{proof}

\begin{remark} {\em For positive $l$ one can use the homomorphism $\delta\colon
H^2_c(k,B)_l\to \Hb^2(B)_l$ to obtain a similar description for
$\bigoplus_{l>0} \Hb^2(B)_l$.
However, the positive part of the truncated bialgebra cohomology is not
relevant in the
context of graded bialgebra deformations.}
\end{remark}

\subsection{Graded bialgebra deformations}
Now let $B=u({\mathcal D},\lambda,\mu)$, defined in Section \ref{sec-fdpha}.
Also assume that the order of $\Gamma$ is not divisible by primes
$<11$. These Hopf algebras are in general filtered by the
coradical filtration, with $\deg(x_i)=1$ ($i=1,\ldots,\theta$) and
$\deg(g)=0$ ($g\in \Gamma$). The filtration allows us to define
related Hopf algebras over $k[t]$, where $t$ is an indeterminate,
as in \cite{DCY}: By \cite[Theorem 3.3(1)]{AS8}, $B$ has
 PBW basis $\{x_{\beta_1}^{a_1}\cdots x_{\beta_p}^{a_p} g\mid
1\leq a_i < N_i, \ g\in\Gamma\}$.
Express each element of $B$ uniquely as a linear combination
of these basis elements.
Then there exist unique maps $m_s:B\ot B\rightarrow B$,
homogeneous of degree $-s$, such that
$$
  m(a\ot b) = \sum_{s\geq 0} m_s(a\ot b)
$$
for all $a,b\in H$.
Now define a new multiplication $m_t:B\ot B\rightarrow B[t]$ by
$$
   m_t(a\ot b) = \sum_{s\geq 0} m_s(a\ot b) t^s,
$$
and extend $k[t]$-linearly to $B[t]\ot_{k[t]}B[t]$. In particular
the analogs of the  linking and root vector relations
(\ref{linking}) and (\ref{rootvector}) for $B[t]$ will now involve
powers of $t$. When we write $B[t]$, we will always mean the
vector space $B[t]$ with multiplication $m_t$ and the usual
(graded) comultiplication. In this way $B[t]$ is a graded
deformation of $\Gr B = u(\mathcal{D},0,0)$.

If $s>0$, then define $\lambda^{(s)}$ and $\mu^{(s)}$ by
$$\lambda^{(s)} =
\begin{cases} \lambda ;&\mbox{ if } s=2 \\ 0 ;&\mbox{ otherwise
}\end{cases},\; \mu^{(s)}_{\alpha} = \begin{cases} \mu_\alpha
;&\mbox{ if } s = \Ht(\alpha)N_\alpha \\ 0 ;&\mbox{ otherwise
}\end{cases}.$$ Note that $u(\mathcal{D},\mu,\lambda)[t]/(t^{s+1})
\simeq u(\mathcal{D},\mu^{(s)},\lambda^{(s)})[t]/(t^{s+1})$.

Let $r=r(\lambda,\mu)$ be the smallest positive integer, if it
exists, such that $m_r\neq 0$. (If it does not exist, set
$r=r(\lambda,\mu)=0$.) Since $B[t]$ is a bialgebra, $(m_r,0)$ is
necessarily a bialgebra two-cocycle, of degree $-r$, and
$B[t]/(t^{r+1})$ is an $r$-deformation (see Section \ref{gbc}).
Note that if $\lambda\neq 0$, then $r=2$, due to the linking
relations (\ref{linking}). Recall the definitions (\ref{lambdaji})
of $\lambda_{ji}$ ($i<j$) and of $F_{ji}, F_{\alpha}$ via
(\ref{F2}) above. The following is a nice consequence of Theorem
\ref{Fbasis}.

\begin{thm}\label{thm-def} \indent\\
\begin{itemize}
\item[(i)] Let $B[t]=u(\mathcal{D},\lambda,\mu)[t]$ and
$B'[t]=u(\mathcal{D},\lambda',\mu')[t]$ and let $s$ be a
nonnegative integer such that $\lambda^{(s)}=\lambda'^{(s)}$ and
$\mu^{(s)}=\mu'^{(s)}$. Then $B[t]/(t^{s+2})$ and
$B'[t]/(t^{s+2})$ are $(s+1)$-deformations of $u(\mathcal{D},0,0)$
extending the same $s$-deformation $B[t]/(t^{s+1}) =
B'[t]/(t^{s+1})$. If $m_t = m + tm_1 +\ldots t^s m_s + t^{s+1}
m_{s+1} +\ldots$ then $m'_t = m + tm_1+\ldots \ldots t^s m_s +
t^{s+1}m'_{s+1}+\ldots$ and $(m_{s+1}-m'_{s+1},0)$ is a bialgebra
cocycle of degree $-(s+1)$ cohomologous to $(F,0)$, where
$$
 F = \sum_{1\leq i<j\leq \theta} \delta_{s+1,2}(\lambda_{ji}-\lambda'_{ji})F_{ji} +
     \sum _{\alpha\in\Phi^+} \delta_{s+1,\Ht(\alpha)N_{\alpha}}
          (\mu_{\alpha}-\mu'_{\alpha}) F_{\alpha}.
$$

\item[(ii)] In particular, the Hopf algebra $u({\mathcal
D},\lambda,\mu)[t]$ is a graded deformation of $u({\mathcal
D},0,0)= R\# k\Gamma$, over $k[t]$, with infinitesimal deformation
$(m_r,0)$ cohomologous to $(F,0)$, where
$$
  F =  \sum_{1\leq i<j\leq \theta} \delta_{r,2}\lambda_{ji}F_{ji} +
     \sum _{\alpha\in\Phi^+} \delta_{r,\Ht(\alpha)N_{\alpha}}
          \mu_{\alpha} F_{\alpha}.
$$
\end{itemize}
\end{thm}
\begin{proof}
In view of Theorem \ref{Fbasis} and Remark \ref{rdiff} it is clear
that $(m_{s+1}-m'_{s+1},0)$ is a bialgebra cocycle cohomologous to
$(F,0)$, where $F=\sum a_{ji} f_{ji} + \sum b_\alpha f_\alpha$ for
some scalars $a_{ji}, b_\alpha$. Evaluating $m_{s+1}$ and
$m'_{s+1}$ at $x_i\ot x_j$ (if $s=1$) and $x_\alpha\ot
x_\alpha^{N_\alpha-1}$ (if $s=\Ht(\alpha)N_\alpha -1$) then
identifies these coefficients.
\end{proof}

We give one more class of examples to which our cohomological
techniques apply, the rank one Hopf algebras of Krop and Radford
\cite{KR}. Assume $k$ has characteristic 0. Let $\theta=1$ and
$(a_{11})=(2)$. Let $\Gamma$ be a finite group (not necessarily
abelian), $a=g_1$ a central element of $\Gamma$, and
$\chi\in\widehat{\Gamma}$. Let $N$ be the order of $\chi(a)$. Let
$x=x_1$, and $R=k[x]/(x^N)$, on which $\Gamma$ acts via $\chi$,
that is $g\cdot x =\chi(g)x$ for all $g\in\Gamma$. Let $B=R\#
k\Gamma$, a generalized Taft algebra with $\Delta(x) =x\ot 1 +
a\ot x$. Similar to the functions $f_{\alpha}$ in Section 4.1,
there is a Hochschild two-cocycle $f:R\ot R\rightarrow k$ defined
by
$$
  f(x^i,x^j) = \left\{ \begin{array}{rl}
                1 & \mbox{ if } i+j=N\\
                0 & \mbox{ otherwise}.\end{array}\right.
$$
This cocycle is $\Gamma$-invariant precisely when $\chi^N=\varepsilon$.
In this case, let $\mu\in k$.
There is a bialgebra deformation of $B$ in which the relation
$x^N=0$ is deformed to $x^N = \mu (1-a^N)$;
this is the Hopf algebra $H_{\mathcal{D}}$ of Krop and Radford
\cite{KR}. In case $\Gamma$ is abelian, this example is included
in the Andruskiewitsch-Schneider classification.


\begin{thebibliography}{9999}

\bibitem{AS5}  N.\ Andruskiewitsch and H.-J.\ Schneider,
``Pointed Hopf algebras,'' in: New Directions in Hopf Algebras,
MSRI Publications 43, 1--68, Cambridge Univ.\ Press, 2002.

\bibitem{AS8}  N.\ Andruskiewitsch and H.-J.\ Schneider,
``On the classification of finite-dimensional pointed Hopf algebras,''
math.QA/0502157, to appear in Ann.\ Math.

\bibitem{BG} A.\ Braverman and D.\ Gaitsgory,
``Poincar\'{e}-Birkhoff-Witt Theorem for quadratic algebras of Koszul
type,''
J.\ Algebra 181 (1996), 315--328.

\bibitem{CGW} A.\ C\u{a}ld\u{a}raru, A.\ Giaquinto, and S.\ Witherspoon,
``Algebraic deformations arising from orbifolds with discrete torsion,''
J.\ Pure Appl.\ Algebra 187 (2004), 51--70.

\bibitem{DCK} C.\ De Concini and V.\ G.\ Kac, ``Representations
of quantum groups at roots of 1,'' in {\em Operator Algebras, Unitary
Representations, Enveloping Algebras, and Invariant Theory: actes du
Colloque en l'honneur de Jacques Dixmier}, ed.\ A.\ Connes et al.,
Progr.\ Math.\ 92, Birkh\"{a}user, Boston, 1990, 471--506.

\bibitem{DCY} Y.\ Du, X.-W.\ Chen, and Y.\ Ye,
``On graded bialgebra deformations,'' Algebra Colloq.\ 14 (2007), no.\ 2,
301--312.

\bibitem{GS} M.\ Gerstenhaber and  S.\ D.\ Shack, ``Bialgebra cohomology,
deformations, quantum
groups and algebraic deformations,'' Proc.\ Natl.\ Acad.\ Sci.\ USA, vol.\
87 (1990), 478--481.

\bibitem{GS2}  M.\ Gerstenhaber and  S.\ D.\ Shack, ``Algebras,
bialgebras,
quantum groups, and algebraic deformations,'' in {\em Deformation theory
and quantum groups with applications to mathematical physics}, 51--92,
Contemp.\ Math.\ 134, Amer.\ Math.\ Soc., 1992.

\bibitem{GK} V.\ Ginzburg and S.\ Kumar, ``Cohomology of quantum groups
at roots of unity,''  Duke Math.\ J.\ 69 (1993), no.\ 1, 179--198.

\bibitem{Gr} L.\ Grunenfelder, ``Tangent cohomology, Hopf algebra actions and
deformations,'' J.\ Pure Appl.\ Algebra 67 (1990), 125–149.

\bibitem{GM} L.\ Grunenfelder and M.\ Mastnak, ``Cohomology of abelian
matched
pairs and the Kac sequence,'' J.\ Algebra 276  (2004),  no. 2, 706--736.

\bibitem{H} I.\ Heckenberger, ``Classification of arithmetic root
systems,''
math.QA/0605795.

\bibitem{kac90} V.\ G.\ Kac, {\em Infinite dimensional Lie algebras},
Cambridge University Press, 3d ed., 1990.

\bibitem{KR} L.\ Krop and D.\ E.\ Radford, ``Finite-dimensional
Hopf algebras of rank one in characteristic zero,'' J.\ Algebra
302 (2006), no.\ 1, 214--230.

\bibitem{L} G.\ Lusztig, ``Quantum groups at roots of 1,''
Geom.\ Dedicata 35 (1990), 89--114.

\bibitem{Mas2} A.\ Masuoka, ``Abelian and non-abelian second cohomologies
of quantized enveloping algebras,'' math.QA/0708.1982.

\bibitem{RA} D.\ E.\ Radford, ``The structure of Hopf algebras with a projection,''
 J.\ Algebra 92 (1985), no.\ 2, 322--347.

\bibitem{Sc}P.\ Schauenburg, ``Hopf modules and Yetter-Drinfeld modules,''
J. Algebra  169  (1994),  no. 3, 874--890.

\bibitem{S} D.\ \c{S}tefan, ``Hochschild cohomology on Hopf Galois
extensions,''
J.\ Pure Appl.\ Algebra 103 (1995), 221--233.

\bibitem{T} R.\ Taillefer, ``Cohomology theories of Hopf bimodules and
cup-product,''
Algebr. Represent. Theory  7  (2004),  no. 5, 471--490.

\bibitem{W} C.\ Weibel, An Introduction to Homological Algebra, Cambridge
University Press
(1994).

\end{thebibliography}
\end{document}